\documentclass[12pt,a4paper,leqno]{article}
\usepackage{etoolbox} \usepackage{comment}
\usepackage[margin=1in]{geometry}
\usepackage{times}
\usepackage[amsthm,thmmarks]{ntheorem}
\usepackage{amsfonts,amssymb}
\usepackage{amsmath,xypic}
\usepackage{pdflscape}
\usepackage{tocloft}
\usepackage[export]{adjustbox} 
\usepackage[cal=boondoxo,calscaled=.96,scr=rsfs]{mathalpha} \usepackage{multicol}
\usepackage{hyphenat}
\usepackage[usenames,dvipsnames]{color}
\usepackage{rotating}
\usepackage{longtable}
\usepackage{caption}
\usepackage{tcolorbox}

\usepackage{epigraph}
\usepackage{enumitem}
\setlist[enumerate]{listparindent=0.5in}
\newcommand{\be}{\begin{equation}}
\newcommand{\ee}{\end{equation}}
\newcommand{\bes}{\begin{equation*}}
\newcommand{\ees}{\end{equation*}}
\newcommand{\bea}{\begin{eqnarray}}
\newcommand{\eea}{\end{eqnarray}}
\newcommand{\beas}{\begin{eqnarray}}
\newcommand{\eeas}{\end{eqnarray}}
\newcommand{\ben}{\begin{note}}
\newcommand{\een}{\end{note}}
\newcommand{\bexl}{\vskip0.1em\noindent\hrulefill\vskip1em\begin{ExerciseList}}
\newcommand{\eexl}{\end{ExerciseList}\hrulefill}

\newcommand{\bthm}{\begin{theorem}}
\newcommand{\ethm}{\end{theorem}}
\newcommand{\bpro}{\begin{prop}}
\newcommand{\epro}{\end{prop}}
\newcommand{\bcor}{\begin{corollary}}
\newcommand{\ecor}{\end{corollary}}
\newcommand{\bcon}{\begin{conjecture}}
\newcommand{\econ}{\end{conjecture}}
\newcommand{\bp}{\begin{proof}}
\newcommand{\ep}{\end{proof}}
\newcommand{\blem}{\begin{lemma}}
\newcommand{\elem}{\end{lemma}}
\newcommand{\bn}{\begin{note}}
\newcommand{\en}{\end{note}}
\newcommand{\benum}{\begin{enumerate}}
\newcommand{\eenum}{\end{enumerate}}
\newcommand{\bed}{\begin{defn}}
\newcommand{\eed}{\end{defn}}
\newcommand{\brem}{\begin{remark}}
\newcommand{\erem}{\end{remark}}

\newcommand{\btik}{\begin{tikzpicture}\begin{axis}[scale=0.5,axis y line=center, axis x line=middle]}
\newcommand{\etik}{\end{axis}\end{tikzpicture}}

\let\into=\hookrightarrow
\let\mapsto=\longmapsto

\newcommand{\upperRomannumeral}[1]{\uppercase\expandafter{\romannumeral#1}}

 \usepackage{tikz}
\usetikzlibrary{mindmap,backgrounds}
\usetikzlibrary{automata}
\usepackage{graphicx}
\usepackage[pagewise]{lineno}
\usepackage{stackengine}
\usepackage{stmaryrd}
\usepackage{microtype}

\usepackage{titlesec}
\usepackage{url}
\usepackage{natbib}
	\let\cite=\citep
\usepackage{fancyhdr}
\usepackage[colorlinks,citecolor=blue]{hyperref}

\newtoggle{draft}
\togglefalse{draft}

\newcommand{\preliminary}{{\\ \normalsize \textcolor{blue}{Preliminary version for comments}}{\relax}}
\newtheorem{theorem}[equation]{Theorem}      \newtheorem{lemma}[equation]{Lemma}          \newtheorem{corollary}[equation]{Corollary}  \newtheorem{proposition}[equation]{Proposition}

\theoremstyle{definition}
\newtheorem{conj}[equation]{Conjecture}

\theoremstyle{definition}
\newtheorem{defn}[equation]{Definition}
\theoremstyle{remark}

\theoremsymbol{\ensuremath{\bullet}}

\theoremstyle{definition}
\newtheorem{remark}[equation]{Remark}
\theoremsymbol{\ensuremath{\bullet}}

\numberwithin{equation}{section}

\newcommand{\para}{\subsection{}\nwss}

\titleformat{\subsection}[runin]{\normalfont\bfseries}{\S\ \thesubsection}{.5em}{}[{\ \ }]
\titlespacing{\subsection}{0pt}{1.5ex plus .1ex minus .2ex}{0pt}
\newcommand{\subpara}{\subsubsection{}\nwsss}

\titleformat{\subsubsection}[runin]{\normalfont\bfseries}{\S\ \thesubsubsection}{.5em}{}[{\ \ }]
\titlespacing{\subsubsection}{0pt}{1.5ex plus .1ex minus .2ex}{0pt}

\let\into=\hookrightarrow
\let\isom=\simeq

\let\tensor=\otimes
\newcommand{\A}{\mathscr{A}}
\newcommand{\abs}[1]{\left\vert#1\right\vert}

\newcommand{\bF}{{\bar{F}}}
\newcommand{\bQ}{{\bar{\Q}}}

\newcommand{\C}{{\mathbb C}}

\newcommand{\F}{{\mathbb F}}

\newcommand{\gal}{{\rm Gal}}

\newcommand{\N}{\mathscr{N}}

\newcommand{\Q}{{\mathbb Q}}
\newcommand{\R}{{\mathbb R}}

\newcommand{\Spec}{{\rm Spec}}

\newcommand{\Z}{{\mathbb Z}}

\renewcommand{\O}{{\mathscr O}}
\renewcommand{\P}{{\mathbb P}}
\renewcommand{\wp}{{\mathfrak p}}

\newcommand{\fm}{{\mathfrak{M}}}
\newcommand{\sM}{\mathscr{M}}

\newcommand{\invlim}{\varprojlim}
\let\fm=\fa
\newcommand{\mapright}[1]{{\xymatrix{{}\ar[r]^{#1}&{}}}}
\newcommand{\mapleft}[1]{{\xymatrix{{}&{}\ar[l]_{#1}}}}

\renewcommand{\bpro}{\begin{proposition}}
	\renewcommand{\epro}{\end{proposition}}
\renewcommand{\bcon}{\begin{conj}}
	\renewcommand{\econ}{\end{conj}}

\title{Construction of  Arithmetic Teichmuller spaces II: Proof of a local prototype of Mochizuki's Corollary~3.12
\preliminary
}
\author{Kirti Joshi}

\newcommand{\Address}{\bigskip\noindent{\footnotesize\textsc{{Math. department, University of Arizona, 617 N Santa Rita, Tucson
		85721-0089, USA.}}\par\nopagebreak 
\noindent\textit{Email:}	\texttt{kirti@math.arizona.edu}}}

\begin{document}
	\maketitle
\epigraphwidth0.55\textwidth
\epigraph{Tamso m\=a jyotir gamaya\footnotemark[1]}{Bṛhadāraṇyaka Upaniṣad \cite{olivelle-upanishads}}\footnote{\textit{From darkness, lead me to light.} [Trans. taken from \cite{olivelle-upanishads}.]}

\lhead{}

\iftoggle{draft}{\pagewiselinenumbers}{\relax}
\newcommand{\act}{\curvearrowright}
\newcommand{\lmp}{{\Pi\act\Ot}}
\newcommand{\lmpi}{{\lmp}_{\int}}
\newcommand{\lmpf}{\lmp_F}
\newcommand{\Om}{\O^{\times\mu}}
\newcommand{\Omf}{\O^{\times\mu}_{\bF}}
\renewcommand{\N}{\mathbb{N}}
\newcommand{\yoga}{Yoga}
\newcommand{\gl}[1]{{\rm GL}(#1)}
\newcommand{\bK}{\overline{K}}
\newcommand{\reptrip}{\rho:G_K\to\gl V}
\newcommand{\reptripp}[1]{\rho\circ\alpha:G_{\ifstrempty{#1}{K}{{#1}}}\to\gl V}
\newcommand{\benumlab}{\begin{enumerate}[label={{\bf(\arabic{*})}}]}
\newcommand{\ord}{\mathop{\rm ord}\nolimits}	
\newcommand{\kcs}{K^\circledast}
\newcommand{\lcs}{L^\circledast}
\renewcommand{\A}{\mathbb{A}}
\newcommand{\bfq}{\bar{\mathbb{F}}_q}
\newcommand{\tripod}{\P^1-\{0,1728,\infty\}}

\newcommand{\vseq}[2]{{#1}_1,\ldots,{#1}_{#2}}
\newcommand{\anab}[4]{\left({#1},\{#3 \}\right)\anabelmap\left({#2},\{#4 \}\right)}

\newcommand{\gln}{{\rm GL}_n}
\newcommand{\glo}[1]{{\rm GL}_1(#1)}
\newcommand{\glt}[1]{{\rm GL_2}(#1)}

\newcommand{\iut}{\cite{mochizuki-iut1, mochizuki-iut2, mochizuki-iut3,mochizuki-iut4}}
\newcommand{\topics}{\cite{mochizuki-topics1,mochizuki-topics2,mochizuki-topics3}}

\newcommand{\linv}{\mathfrak{L}}
\newcommand{\bedef}{\begin{defn}}
\newcommand{\eedef}{\end{defn}}
\renewcommand{\act}[1][]{\overset{#1}{\curvearrowright}}
\newcommand{\bfx}{\overline{F(X)}}
\newcommand{\anabelmap}{\leftrightsquigarrow}
\newcommand{\ban}[1][G]{\mathscr{B}({#1})}
\newcommand{\pit}{\Pi^{temp}}
 
 \newcommand{\bL}{\overline{L}}
 \newcommand{\bkm}{\bK_M}
 \newcommand{\vbk}{v_{\bK}}
 \newcommand{\vbkm}{v_{\bkm}}
\newcommand{\ocs}{\O^\circledast}
\newcommand{\ot}{\O^\triangleright}
\newcommand{\ocsk}{\ocs_K}
\newcommand{\otk}{\ot_K}
\newcommand{\ok}{\O_K}
\newcommand{\oko}{\O_K^1}
\newcommand{\oks}{\ok^*}
\newcommand{\Qpb}{\overline{\Q}_p}
\newcommand{\Qpbh}{\widehat{\overline{\Q}}_p}
\newcommand{\tr}{\triangleright}
\newcommand{\ocpt}{\O_{\C_p}^\tr}
\newcommand{\ocpf}{\O_{\C_p}^\flat}
\newcommand{\sG}{\mathscr{G}}
\newcommand{\sxfe}{\mathscr{X}_{F,E}}
\newcommand{\sxfep}{\mathscr{X}_{F,E'}}
\newcommand{\loglt}{\log_{\sG}}
\newcommand{\fc}{\mathfrak{t}}
\newcommand{\ku}{K_u}
\newcommand{\kup}{\ku'}
\newcommand{\kt}{\tilde{K}}
\newcommand{\sGpf}{\sG(\O_K)^{pf}}
\newcommand{\hgm}{\widehat{\mathbb{G}}_m}
\newcommand{\bE}{\overline{E}}
\newcommand{\bEp}{\overline{E}'}
\newcommand{\sY}{\mathscr{Y}}
\newcommand{\syfe}{\mathscr{Y}_{F,E}}
\newcommand{\syfqp}[1]{\mathscr{Y}_{\cptl{#1},\Q_p}}
\newcommand{\syfQp}{\mathscr{Y}_{{F},\Q_p}}
\newcommand{\syfqpe}[1]{\mathscr{Y}_{{#1},E}}
\newcommand{\fJ}{\mathfrak{J}}
\newcommand{\fM}{\mathfrak{M}}
\newcommand{\locvar}{local arithmetic-geometric anabelian variation of fundamental group of $X/E$ at $\wp$}
\newcommand{\fjxep}{\fJ(X,E,\wp)}
\newcommand{\fjxe}{\fJ(X,E)}
\newcommand{\fpc}[1]{\widehat{{\overline{\F_p(({#1}))}}}}
\newcommand{\cpt}{\C_p^\flat}
\newcommand{\cptl}[1]{\C_{p,{#1}}^\flat}
\newcommand{\fja}[1]{\fJ^{\rm arith}({#1})}
\newcommand{\ainfe}{A_{\inf,E}(\O_F)}
\renewcommand{\ainfe}{W_{\O_E}(\O_F)}
\newcommand{\gmh}{\widehat{\mathbb{G}}_m}
\newcommand{\sE}{\mathscr{E}}
\newcommand{\bpi}{B^{\varphi=\pi}}
\newcommand{\bpip}{B^{\varphi=p}}
\newcommand{\onto}{\twoheadrightarrow}

\newcommand{\cpmax}{{\C_p^{\rm max}}}
\newcommand{\xan}{X^{an}}
\newcommand{\yan}{Y^{an}}
\newcommand{\bPi}{\overline{\Pi}}
\newcommand{\bPit}{\bPi^{\rm{\scriptscriptstyle temp}}}
\newcommand{\Pit}{\Pi^{\rm{\scriptscriptstyle temp}}}
\renewcommand{\pit}[1]{\Pi^{\scriptscriptstyle temp}_{#1}}
\newcommand{\pitb}[1]{{\overline{\Pi}}^{\scriptscriptstyle temp}_{#1}}
\newcommand{\pitk}[2]{\Pi^{\scriptscriptstyle temp}_{#1;#2}}
\newcommand{\pio}[1]{\pi_1({#1})}
\newcommand{\fTeich}{\widetilde{\fJ(X/L)}}
\newcommand{\ssep}{\S\,} \newcommand{\vphi}{\varphi}
\newcommand{\sgt}{\widetilde{\sG}}

\newcommand{\flog}{\mathfrak{log}}

\togglefalse{draft}
\newcommand{\FF}{\cite{fargues-fontaine}}
\iftoggle{draft}{\pagewiselinenumbers}{\relax}

\newcommand{\attportion}{Sections~\ref{se:number-field-case}, \ref{se:construct-att}, \ref{se:relation-to-iut}, \ref{se:self-similarity} and \ref{se:applications-elliptic}}

\newcommand{\four}{Sections~\ref{se:grothendieck-conj}, \ref{se:untilts-of-Pi}, and \ref{se:riemann-surfaces}}

\numberwithin{equation}{subsection}

\newcommand{\tcp}{\widetilde{\C}_p}
\newcommand{\tK}{\widetilde{K}}
\newcommand{\tF}{\widetilde{F}}
\newcommand{\bfp}{\overline{\F}_p}
\newcommand{\sxqp}{\mathscr{X}_{\cpt,\Q_p}}
\newcommand{\syQp}{\mathscr{Y}_{\cpt,\Q_p}}
\newcommand{\sxQp}{\mathscr{X}_{\cpt,\Q_p}}
\newcommand{\ttxlv}{\tilde{\Theta}_{X,\ell;v}}
\newcommand{\tm}[1]{\theta_M({#1})}
\renewcommand{\fm}{\frak{m}}
\newcommand{\ells}{{\ell^*}}
\newcommand{\cpmaxt}{\C_p^{max\flat}}
\newcommand{\sL}{\mathscr{L}}
\newcommand{\sigfjxe}{\Sigma\fjxe_{F}}
\newcommand{\sigfjxecpt}{\Sigma\fjxe_{\cpt}}

\newcommand{\moccor}{\cite[Corollary 3.12]{mochizuki-iut3}}
\newcommand{\nwss}{\numberwithin{equation}{subsection}}
\newcommand{\nwsss}{\numberwithin{equation}{subsubsection}}

\nwss

\newcommand{\ttheta}{\widetilde{\Theta}}

\begin{abstract}
		This paper deals with  consequences of the existence of Arithmetic Teichmuller spaces established \href{https://arxiv.org/abs/2106.11452}{here} and \href{https://arxiv.org/abs/2010.05748}{here}. Theorem~\ref{th:main} provides a proof of a local version of  Mochizuki's Corollary~3.12. Local means for a fixed $p$-adic field. There are several new innovations in this paper. Some of the main results are as follows. Theorem~\ref{th:tate-parameter-as-a-function} shows that one can view the Tate parameter of Tate elliptic curve as a function on the arithmetic Teichmuller space of \cite{joshi-teich}, \cite{joshi-untilts}. The next important point is the  construction of Mochizuki's $\Theta_{gau}$-links  and the set of such links, called Mochizuki's Ansatz in \ssep\ref{se:ansatz}. Theorem~\ref{pr:lift-vals} establishes  valuation scaling property satisfied by points of Mochizuki's Ansatz (i.e. by my version of $\Theta_{gau}$-links). These results lead to the construction of a theta-values set (\ssep\ref{se:construction-ttheta}) which is similar to Mochizuki's Theta-values set (differences between the two are in \ssep\ref{pa:diff-between-two-sets}). Finally Theorem~\ref{th:main} is established. For completeness, I provide an intrinsic proof  of the existence of Mochizuki's $\log$-links (Theorem~\ref{th:log-link}),  $\flog$-links (Theorem~\ref{th:flog-links}) and Mochizuki's log-Kummer Indeterminacy (Theorem~\ref{th:flog-kummer-correspondence}) in my theory.
\end{abstract}

\setcounter{tocdepth}{1}

\tableofcontents

\section{Introduction} \para  In \cite{joshi-teich}, \cite{joshi-untilts}, I have detailed a theory of arithmetic Teichmuller spaces arising from arithmetic and providing an arithmetic analog of classical holomorphic structures. While \cite{joshi-teich} details a Berkovich function theoretic approach to arithmetic holomorphic structures,   \cite{joshi-untilts} details the group theoretic aspect of my theory which allows one to apply my results to Mochizuki's group theoretic viewpoint detailed in \iut.

\subpara The important remarkable innovation of \iut\ is its original strategy of bounding arithmetic degrees of interest in the proofs of Vojta and Szpiro conjectures in \cite{mochizuki-iut4}. The key result which underpins this strategy is \moccor.  The center piece of \moccor\ (and a consequence of \cite[Theorem 3.11]{mochizuki-iut3}) is  the construction of a certain set which I will denote by $\ttheta_{Mochizuki}$ here.   The assertion of \moccor\ provides a lower bound on a suitably defined notion of the size of $\ttheta_{Mochizuki}$.  

\subpara The said corollary should be understood as an averaging procedure for functions on a suitable configuration space of arithmetic holomorphic structures at each prime for a fixed hyperbolic curve over a number field. The space of arithmetic holomorphic structures is the Arithmetic Teichmuller Space constructed in \cite{joshi-teich,joshi-untilts}. This has many properties similar to classical Teichmuller spaces which appear in the adelic theory as archimedean contributions. \textit{Notably this space determines a unique isomorphism class of a profinite  group given as the \'etale fundamental group of the each point  of this space} (for classical Teichmuller spaces this corresponds to $\Sigma'\in T_{\Sigma}\mapsto \pi_1(\Sigma')$). This is the perspective of the present series of papers. Similar averaging procedures exist in classical Teichmuller Theory  and are quite well-understood. The  archimedean case of the said corollary is discussed in \cite{joshi-teich-quest} which also provides a generally accessible  discussion of \moccor.

\subpara Important consequence of the existence of arithmetic Teichmuller spaces of \cite{joshi-teich,joshi-untilts} is that the Tate parameter (and hence values of theta-functions at suitably chosen points) of an elliptic curve with potentially multiplicative reduction may be considered as a function  on its arithmetic Teichmuller space (Theorem~\ref{th:tate-parameter-as-a-function}) and \moccor\ should be viewed as an averaging of this function over a suitable configuration space of the arithmetic Teichmuller space. This function (in the sense of algebraic geometry) does not take value in a common valued field. This necessitates a lifting procedure (detailed in \ssep\ref{se:lifting-values-to-B}) which allows one to consider lifts of the values in a fixed ring provided by \cite{fargues-fontaine}.

\subpara Mochizuki refers to his set $\ttheta_{Mochizuki}$ constructed in \moccor\ as a \textit{multi-radial representation of theta-values}. This term  simply means that $\ttheta_{Mochizuki}$ is  the set of theta-values (computed in a fixed chosen set) arising from distinct arithmetic holomorphic structures. Hence the construction of $\ttheta_{Mochizuki}$ is predicated on the existence of distinct arithmetic holomorphic structures. Additionally the construction of $\ttheta_{Mochizuki}$ is also predicated on the existence of $\Theta_{gau}$-Links discovered by Mochizuki \cite{mochizuki-iut2}. 

\subpara  In \cite{joshi-teich}, \cite{joshi-untilts} I have demonstrated the existence of distinct arithmetic holomorphic structures and as detailed in \cite{joshi-untilts}, these also provide distinct arithmetic holomorphic structures in the sense of \iut. 

An important discovery of this paper is the intrinsic geometric construction of $\Theta_{gau}$-Links  and the set of such links (I do not refer to this construction as $\Theta_{gau}$-Links). This is detailed in \ssep~\ref{se:ansatz} with the principal property of $\Theta_{gau}$-Links established in Theorem~\ref{pr:lift-vals}.

\subpara This paper details my approach to the construction of sets of type $\ttheta_{Mochizuki}$. In \ssep\ref{se:construction-ttheta}, I construct the set $\ttheta$ which is corresponds to $\ttheta_{Mochizuki}$ in my theory. The difference between the two sets is discussed in \ssep\ref{pa:diff-between-two-sets}.  One of the main results of this paper is Theorem~\ref{th:main} which  establishes a lower bound for the set $\ttheta$  similar to \moccor. In particular Theorem~\ref{th:main} is a prototype (at one prime) of the bound established in \moccor.  Owing to the difficult presentation of the claims and results of \iut,  the existence of $\ttheta_{Mochizuki}$ which is asserted in \moccor\ has proved quite difficult to verify independently and even its existence has  evoked a sense of incredulity and  disbelief among many readers of \iut. 

\subpara Hence my  main purpose in establishing the prototype result (Theorem~\ref{th:main}) is to demonstrate that sets ($\ttheta$ and $\ttheta_{Mochizuki}$)  and the bounds claimed therein  can be indeed contemplated and rigorously established once one has the results of the present paper, the framework of arithmetic Teichmuller spaces detailed in this paper and \cite{joshi-teich}, \cite{joshi-untilts}.  As is pointed out in \cite{joshi-untilts} combining my results with Mochizuki's group theoretic framework, and incorporating the changes to \cite{mochizuki-iut3} discussed in \cite{joshi-untilts},  one may also obtain the corresponding claims using Mochizuki's group theoretic framework. The general case of this corollary is treated in \cite{joshi-teich-rosetta}.

\subpara  One may also construct $\ttheta_{Mochizuki}$ using methods of this paper and \cite{joshi-teich,joshi-untilts}. In particular at this juncture there is no doubt that multi-radial representations of theta-values as asserted in \cite{mochizuki-iut3} exists, even if one works with number fields, as the existence is  separately established for each prime (and this is carried out here).

\para Let me now discuss  Mochizuki's \textit{$\Theta_{gau}$-Link} and my proof of its existence.  
\newcommand{\tsigcpt}{\tilde{\Sigma}_{\cpt}}
\newcommand{\tsigf}{\tilde{\Sigma}_{F}}
\newcommand{\sigcpt}{{\Sigma}_{\cpt}}
\newcommand{\sigf}{{\Sigma}_{F}}

\subpara Mochizuki's $\Theta_{gau}$-Link rests on the remarkable assertion that  if $q$ is the Tate parameter of a Tate elliptic curve $X/E$ over a $p$-adic field, then  one can set (intentionally written in quotes--for this discussion) $$\text{``}q=q^{j^2}\text{''} \text{ and } \abs{q}_{\Q_p}<1$$ for integers $j$ and even simultaneously set
\be\label{eq:ansatz} \text{``}q=(q,q^{2^2},\cdots, q^{\ells^2}).\text{''}
\ee 
Mochizuki has asserted  presciently (see \cite[Pages 408--409 and the table on Page 411]{mochizuki-iut3}), but without proof, that $\Theta_{gau}$-Links (and variants) should be considered as arising from a  Witt-vector type deformation i.e. passage from  $\bmod({p^n})$ to $\bmod({p^{n+1}})$ (for all $n$).  Second important point about $\Theta_{gau}$-Links which is not stated in \iut\ but used there, is that \eqref{eq:ansatz} (and hence the $\Theta_{gau}$-Link) must be understood as providing a simultaneous (but non-trivial) scaling operation in the theory. 

\subpara A key discovery of this paper is that Mochizuki's ansatz \eqref{eq:ansatz}  has an elegant  formulation (\ssep\ref{se:ansatz}):
\be\label{eq:ansatz2} 
 p=[p_K^\flat],p=[p_K^\flat]^{2^2},\cdots,p=[p_K^\flat]^{\ells^2},
\ee
or more succinctly, in the style of \eqref{eq:ansatz} as $$p=([p_K^{\flat}],[(p_K^{\flat})^{2^2}],\cdots,[(p_K^{\flat})^{\ells^2}])$$
where $([p_K^\flat]-p)\subset W(\O_{\cpt})$ is the prime ideal giving rise to some characteristic zero untilt $(K,K^\flat\isom\cpt)$ of $\cpt$. Here  $\cpt$ is the tilt of $\C_p$--see \cite[Lemma 3.4]{scholze12-perfectoid-ihes}--note that $\cpt$ is an algebraically closed perfectoid (hence a complete valued) field of characteristic $p>0$. The system of equations \eqref{eq:ansatz2} define an $\ells$-tuple of closed classical points of $\syQp$. This leads to the construction of the subset $\tsigcpt\subset\abs{\syQp}^\ells$. Each point of $\tsigcpt$ is defined by a local equation to the type \eqref{eq:ansatz2}. The set $\tsigcpt$ will be referred to as \textit{Mochizuki's Primitive Ansatz in $\syQp^\ells$} (I will often conflate $\abs{\syQp}$ with $\syQp$ for notational simplicity)  and it leads to a set of tuples of arithmetic holomorphic structures $\sigfjxe\subset\fjxe^\ells$ (called \textit{Mochizuki's Ansatz}) lying over points of $\tsigcpt$ and which may be understood as the set of $\Theta_{gau}$-Links in the theory of the present paper (see \ssep\ref{se:ansatz}). This allows me to prove  precise versions of both of Mochizuki's  claims regarding $\Theta_{gau}$-links (see Theorem~\ref{pr:lift-vals}).

\subpara By Theorem~\ref{th:tate-parameter-as-a-function} the Tate parameter and hence values of theta functions \ssep\ref{se:theta-values} do not lie in a fixed set. So the problem of comparing them is ill-posed. My approach to this resolving this issue is to consider lifts of the theta values to the ring $B$. In \ssep\ref{se:lifting-values-to-B} I detail my approach to lifting theta-values arising from distinct arithmetic holomorphic structures to the ring $B$. This was first sketched in \cite{joshi-teich}. This allows one to view all the (lifts of)  theta values in one common ring.

\subpara Once the lifting formalism is in place, one can apply it to lifting theta-values provided by a point of Mochizuki's Ansatz to the free $B$-module $B^\ells$. This is the content of Theorem~\ref{thm:theta-pilot-object-appears} which provides, for each point of $\sigfjxecpt$, a collection of tuples
\be\label{eq:theta-pilot-intro} 
([x_1]+\tau_1, \ldots, [x_\ells]+\tau_{\ells})\in B^\ells.
\ee
(see Theorem~\ref{thm:theta-pilot-object-appears} for a precise definition of such tuples; the ring $B=B_{\cpt,\Q_p}$ is the fundamental ring of \cite{fargues-fontaine}). 

Each tuple \eqref{eq:theta-pilot-intro} is the additive analog in my theory (at one prime of a number field) of Mochizuki's multiplicative $\Theta$-Pilot object in \cite{mochizuki-iut3}.  

\subpara The set  $\ttheta$ is the set of tuples \eqref{eq:theta-pilot-intro} for each point  of $\sigfjxecpt$. In other words, by my construction, the set $\ttheta$ is the set of all $\Theta$-Pilot objects i.e. of tuples of theta-values arising from distinct tuples of arithmetic holomorphic structures given by points of $\sigfjxecpt$.  \textit{The set $\ttheta$ is the analog (at one prime) in my theory of Mochizuki's multi-radial representation of theta-values in \moccor}. In particular multi-radial representations of theta-values exists! This is also asserted in \cite{mochizuki-iut3} via  \cite[Theorem 3.11]{mochizuki-iut3} which provides the existence of multi-radiality i.e. of many arithmetic holomorphic structures subject to indeterminacies and \moccor\ which defines Mochizuki's set $\ttheta_{Mochizuki}$ of multi-radial theta-values and provides a lower bound on its size.

\subpara The central difference between the two constructions  ($\ttheta$ and $\ttheta_{Mochizuki}$) is this: Mochizuki uses tuples of Galois cohomology classes in $H^1(G_E,\Q_p(1))^\ells$ to capture his theta-cohomology-classes. By \cite[Proposition 1.4]{mochizuki-theta}, Mochizuki's Galois theoretic theta classes are generated by the values of a theta function which is chosen in  \ssep\ref{se:theta-values} to be in agreement with Mochizuki's choice. This ensures that the  chosen values of this theta function coincide in both the theories.  In the present paper I work with the lifts of these theta-values (i.e. with tuples \eqref{eq:theta-pilot-intro}) to the fundamental ring $B$ of $p$-adic Hodge Theory rather than the Galois cohomology classes they generate.  
From my point of view (and that of \cite{fargues-fontaine}) the ring $B$ (and its variants such as $B_E$) equipped with its canonical action of Galois group  may be considered as the primogenitor of Galois cohomology $H^1(G_{E},\Q_p(1))$ and hence one may apply my results (and methods) to  \iut. Moreover working with the ring $B$ (and its variants) provides a conceptually and algebraically natural way to intrinsically understand the construction of such multi-radial representations.  My results can also be applied to \iut\ using \cite{joshi-untilts}.

\para  An important technical device which plays a central role in \iut, is the notion of \emph{$\log$-links} and a special case of this  called $\flog$-links in \cite{mochizuki-iut3}. Let me now explain how my Arithmetic Teichmuller Theory provides these devices in an intrinsic way. 

\subpara In \ssep\ref{se:log-links}, I prove the existence of a version of  Mochizuki's $\log$-links (Theorem~\ref{th:log-link}) and in Theorem~\ref{th:flog-links} I prove the existence of $\flog$-links from the point of view of  Arithmetic Teichmuller theory of \cite{joshi-teich} \cite{joshi-untilts}. This allows me to demonstrate the precise way in which Mochizuki's $\log$-Kummer correspondence arises (Theorem~\ref{th:flog-links}). 

\subpara Mochizuki has asserted \cite[Page 408]{mochizuki-iut3} that $\flog$-link must be considered to be analogous to the Frobenius morphism in characteristic $p>0$. I demonstrate (Theorem~\ref{th:flog-links}) that this  stated analogy is in fact true quite literally. Especially, I show that a vertical column of $\flog$-links in Mochizuki's $\flog{\rm-}\Theta$-lattice can be identified with a collection of arithmetic holomorphic structures (in the sense of \ssep\ref{se:tate-setup} or \cite{joshi-untilts})
\be\label{eq:log-link-and-canonical-fiber}\left\{(X/E,(E\into K_{y_n}, K_{y_n}^\flat\isom \cpt), *_{K_{y_n}}:\sM(K_{y_n})\to \xan_{E}))\right\}_{n\in\Z}\ee  where  $\{y_n \}_{n\in\Z}\subset\syQp$ is the fiber of the  morphism $\syQp\to\sxQp$ over the canonical closed classical point of $x_{can}\in\sxQp$ (for the existence of this point see \cite[Proposition 10.1.1]{fargues-fontaine} where this point is denoted by $\infty$).

Notably I also prove that $p$-adic valuations of $p$ (and hence of the Tate parameter of $X$) grow in such a collection \eqref{eq:log-link-and-canonical-fiber} of arithmetic holomorphic structures i.e. in a  vertical columns of $\flog$-links. This growth property is analogous to the growth of radii of convergence of solutions of $p$-adic differential equations under Frobenius morphism. This growth property has been suggested by Mochizuki as an analogy (see \cite[Remark 3.6.2, Page 110]{mochizuki-topics3} and compare it with Theorem~\ref{th:flog-kummer-correspondence}\textbf{(4)}).  This  allows me to assert, with a complete degree of certainty, that one has arrived at all the principal landmarks of \cite{mochizuki-iut1,mochizuki-iut2,mochizuki-iut3}.

\subpara Finally the stage is set to prove the prototype of \moccor\ which is a bound for the (suitably defined) size of the set $\ttheta$. This I do in Theorem~\ref{th:main}. 

\subpara The construction of $\ttheta$ carried out here (especially for Theorem~\ref{th:main}) is for a fixed $p$-adic field. Working with elliptic or a once-punctured elliptic curve over a number field requires working simultaneously with finitely many odd primes at which the curve acquires split multiplicative reduction. This requires considerable additional footwork and added notational complexity. But the general version,  adapted for \moccor\ and arithmetic applications is expected to proceed by interpolating between the principal  constructions of this paper, \cite{joshi-teich,joshi-untilts} and \cite{mochizuki-iut3} and is treated in \cite{joshi-teich-rosetta}.

\para At this juncture let me say that by adopting the definitions of arithmetic holomorphic structures given in \cite{joshi-untilts} and the changes this necessitates in some of the portions of \cite{mochizuki-iut3} and using the central discoveries of this paper, one can (for sure) demonstrate the existence of Mochizuki's multi-radial representation of theta values set $\ttheta_{Mochizuki}$.  Moreover it seems \textit{perfectly reasonable} to expect that one can hope to prove \moccor\ (Mochizuki's version) by simply following Mochizuki's method and (his formalism) verbatim with the changes adopted from \cite{joshi-teich}, \cite{joshi-untilts} and the present paper. To facilitate this, in \cite{joshi-teich-rosetta}, I provide a `Rosetta Stone' for translating between geometric objects of my theory and group/Frobenioid theoretic objects of Mochizuki's Theory. On the other hand I have chosen to pursue a different strategy to \moccor\ in my papers because I believe that my geometric approach should have wider applicability as it also works for higher genus, so that one can contemplate higher genus versions of Szpiro's inequality,  and opens the possibility of higher dimensional Diophantine applications.

\subpara \textit{However let me make it absolutely clear that presently I do not make any claims regarding \cite{mochizuki-iut4} simply because personally I have not verified its proofs. My journey in this quest to understand claims of \iut\ has been full of unexpected surprises and hence I do not make any pronouncements on \cite{mochizuki-iut4} at this time.}

\para   For an introduction to the set constructed by Mochizuki's \moccor\  and  other discussions of \moccor\  see \cite{fucheng}, \cite{yamashita},  \cite{dupuy2020statement}. Note that in \cite{dupuy2020statement}, Mochizuki's Corollary 3.12 is stated as Conjecture 1.0.1.

\para  In \cite{mochizuki-iut4}, Mochizuki shows that for his (global  i.e. adelic) version of  the set $\ttheta$ there is an upper bound of the form (my notation here is multiplicative while Mochizuki's notation is additive using the standard (real-valued) logarithm) 
\be\label{eq:mochizukis-inequality} 
\prod_v \abs{q_v^{1/2\ell}}^\ells_{\C_p}\leq { |\ttheta|_{adelic}} \leq A(\ell)\left(\prod_v \abs{q_v^{1/2\ell}}_{\C_p}\right)^{B(\ell)}
\ee where the product is now over all primes (at primes of good reduction one has $q_v=1$), and calculates  $A(\ell)$ and  $B(\ell)$ explicitly showing that $A(\ell)$ depends on $\ell$, the conductor of the underlying elliptic curve and its field of definition while $B(\ell)\geq \ells$ is independent of the curve but depends on $\ell$)  and that this explicit inequality \eqref{eq:mochizukis-inequality} leads to Diophantine applications such as the Szpiro's conjecture (see \iut, \cite{mochizuki-gaussian}, \cite{fesenko-iut},  \cite{fucheng}, \cite{yamashita}, \cite{dupuy2020probabilistic}).

\para {\bf Acknowledgments} 
This paper is a substantially revised version of \cite{joshi-teich-estimates-old-ver} and this acknowledgment is carried over from that \href{https://arxiv.org/pdf/2111.04890.pdf}{version} (one mathematician has requested that he not be thanked for this version of the paper). I am indebted to Taylor Dupuy for reading and providing some comments on early drafts of \cite{joshi-teich-estimates-old-ver} manuscript and for  conversations regarding \iut. I would also like to thank Kiran Kedlaya for some correspondence in the early stage of this work and for his encouragement. Another mathematician who helped me during the course of the work presented here has chosen to remain anonymous, but I would like to acknowledge his help during  the course of this work here.  As I have already stated in \cite{joshi-teich}, neither that paper nor this one, could have existed without Mochizuki's work in anabelian geometry, $p$-adic Teichmuller Theory and his work on the Arithmetic Szpiro Conjecture. So my intellectual debt to Shinichi Mochizuki should be obvious.

\section{Fargues-Fontaine curves $\syQp$ and $\sxQp$ and the ring $B$}\label{se:ff-rings}
\para I will write $F=\cpt$, $\O_F\subset F$ be its valuation ring and $\fm_F\subset \O_F$ be its maximal ideal. This will be fixed throughout this note. Let $\sG/\Z_p$ be a Lubin-Tate formal group with logarithm $\sum_{n=0}^\infty \frac{T^{p^n}}{p^n}$. Then $\sG(\O_F)$ is naturally a Banach space over $\Q_p$ \cite[Chap. 4]{fargues-fontaine}.   Write $B=B_{\Q_p}$ for the ring constructed in \cite[Chap 2]{fargues-fontaine} for the datum $(F=\cpt,E=\Q_p)$.

\para Let $\sxQp$ be the complete Fargues-Fontaine curve, and  $\syQp$ be the incomplete Fargues-Fontaine curve \cite[Chapter 5]{fargues-fontaine}. Of interest to us are the residue fields of closed points $x\in\sxQp$ (resp. $y\in\syQp$). One has a canonical bijection of closed points of degree one \cite{fargues-fontaine}:
$$\abs{\syQp}/\vphi^\Z\to \abs{\sxQp}.$$
If $y\in \abs{\syQp}$ maps to $x\in\sxQp$ under this bijection then one has a natural isomorphism of residue fields
$$K_y\mapright{\isom} K_x.$$

By \cite[Section 10.1]{fargues-fontaine} one has action of 
$\gal(\bQ_p/\Q_p)$ (and hence an action of any of its open subgroups) on $\abs{\syQp}$ (resp. on $\sxQp$). 

By \cite[Theorem 8.29.1]{joshi-teich} one  also has a natural action of ${\rm Aut}_{\Z_p}(\sG(\O_F))$  (resp. ${\rm Aut}_{\Q_p}(\sG(\O_F))$) arising from the natural identification ( \cite[Section 2.3]{fargues-fontaine} (resp. \cite[Th\'eor\`eme 5.2.7]{fargues-fontaine})) of 
$$\abs{\syQp} \isom \left(\sG(\O_F)-\{0\}\right)/\Z_p^*\qquad \text{ resp. } \qquad \abs{\sxQp}\isom \left(\sG(\O_F)-\{0\}\right)/\Q_p^*.$$

\para 
For simplicity I will write $V=\sG(\O_F)$ considered as a topological $\Z_p$-module  and also as a $\Q_p$-Banach space. Both are equipped with an action of $\gal(\bQ_p/\Q_p)$ and in this notation one has 
$$\abs{\sxQp}\isom \P(V), \text{ and}$$
$$\abs{\syQp}\isom (V-\{0\})/\Z_p^*.$$

In this notation the automorphism groups are $${\rm Aut}_{\Z_p}(\sG(\O_F))={\rm Aut}_{\Z_p}(V)$$ and $${\rm Aut}_{\Q_p}(\sG(\O_F))={\rm GL}(V).$$

\para 
The $\gal(\bQ_p/\Q_p)$ action on $\sxQp$ has one unique fixed point  \cite[Proposition 10.1.1]{fargues-fontaine}, denoted by $\infty\in\abs{\sxQp}$ whose residue field is identified with $\C_p$ (which is given the normalized valuation $v_{\C_p}(p)=1$). Similarly one has a canonical point $(\infty,\vphi=1)\in\abs{\syQp}$ whose image in $\sxQp$ is $\infty\in\sxQp$.

So one, by \cite[Theorem 8.29.1]{joshi-teich} and \cite{fargues-fontaine}, has two distinct actions 
$$\gal(\bQ_p/\Q_p)\act \abs{\sxQp} \curvearrowleft {\rm Aut}_{\Q_p}(\sG(\O_F))$$
which move points on the curve and hence their residue fields.

\para As was noted in \cite[\ssep 9.3]{joshi-teich}, the theory of \iut\ is naturally a multiplicative theory. In the multiplicative context of \iut\ one can use the isomorphism
$$\sG(\O_F)\isom \hgm(\O_F).$$

Notably one has $${\rm Aut}_{\Q_p}(\sG(\O_F)) \isom {\rm Aut}_{\Q_p}(\hgm(\O_F))$$
and especially one has the multiplicative action:
$$\gal(\bQ_p/\Q_p)\act \abs{\sxQp} \curvearrowleft {\rm Aut}_{\Q_p}(\hgm(\O_F)).$$

The similarities and differences between this and \iut\ is discussed in \cite[Section 9]{joshi-teich}.

\para Let $F$ be a perfectoid field of characteristic $p>0$. An untilt of $F$ is a pair $(K,\iota:K^\flat\isom F)$ consisting of a perfectoid field $K$ and an isometry $\iota:K^\flat\isom F$ between the tilt $K^\flat$ (\cite{scholze12-perfectoid-ihes}) and $F$. \textit{In this paper the term `untilt of $F$' shall always refer to untilts $(K,\iota:K^\flat\isom F)$ of $F$ in which the fields $K$ have characteristic zero.}
\para\label{pa:value-group-comp-loc} Let me begin by noting that elements of value groups of all untilts of a fixed algebraically closed perfectoid field $F$ of characteristic $p>0$ can be compared in the value group of the tilt $F$. A typical example of $F$ in this paper   is given by $F=\cpt$ which is obtained as the tilt of the completion $\C_p$ of an algebraic closure $\bQ_p$ of $\Q_p$ (\cite{scholze12-perfectoid-ihes}). 
\bpro 
Let $K_1$, $K_2$ be untilts of $\cpt$. . 
\benumlab
\item Then one has $\abs{K_1^*}=\abs{{\C_p^{\flat*}}}$. 
\item Hence, as one has $K_2^\flat=\cpt$,  so $K_1$ and $K_2$ have the same value group.
\item In particular, one can compare elements of $\abs{K_1^*}$ and $\abs{K_2^*}$ as elements of $\abs{{\C_p^{\flat*}}}$. 
\eenum
\epro 
\bp 
The first assertion is  \cite[Lemma 1.3.3]{kedlaya15} or \cite[Lemma 3.4(iii)]{scholze12-perfectoid-ihes}: if $K_1^\flat=\cpt$ then one has $$\abs{K_1^*}=\abs{K_1^{\flat*}}=\abs{{\C_p^{\flat*}}}.$$ 
So if $K_1^\flat=K_2^\flat$ then $\abs{K_1^*}=\abs{K_2^*}=\abs{{\cpt}^*}$. So one can compare elements of value groups of $K_1$ and $K_2$ as elements of the value group of $\cpt$. This proves the proposition.
\ep

\section{Arithmetic Teichmuller Space of a Tate elliptic curve}\label{se:tate-setup}
\para  Let $E$ be a $p$-adic field , $G_E$ be its absolute Galois group for a given algebraic closure of $E$.
Fix a geometrically connected, smooth (quasi-projective) curve $X/E$ where $E$ is a $p$-adic field. Then the objects of the Arithmetic Teichmuller space $\fjxe$ (resp. $\fjxe_{\cpt}$) constructed in \cite{joshi-teich,joshi-untilts} are
$$(Y/E',(E'\into K, K^\flat\isom F), *_{K}:\sM(K)\to \yan_{E'})$$
where  
\benumlab
\item $E'$ is a $p$-adic field, $F$ is some perfectoid field of characteristic $p>0$ (resp. $F=\cpt$), $Y$ is a geometrically connected, smooth (quasi-projective) curve over a $p$-adic field $E'$\item  an untilt $(E'\into K,K^\flat\isom F)$ (resp. $F=\cpt$) of $\cpt$ equipped with an isometric embedding of $E'\into K$, and 
\item a geometric point $\sM(K)\to \yan_{E'}$, and 
\item such that  $Y/E'$ and $X/E$ are tempered anabelomorphic i.e. there exists an isomorphism of tempered fundamental groups $\pi_1^{temp}(Y/E',\sM(K)\to\yan_{E'})\isom \pit{X/E}$ (the fundamental group of $X/E$ computed with some choice of geometric base-point).
\eenum

\para After the results of \cite{joshi-teich} and especially \cite{joshi-untilts} the following definition is reasonable.
\begin{defn} I will refer to the data $(Y/E',(E'\into K, K^\flat\isom F), *_{K}:\sM(K)\to \yan_{E'})$ as an \textit{arithmetic holomorphic structure on $Y/E'$} or simply  as an \textit{arithmetic holomorphoid of $Y/E'$} (or simply as a \textit{holomorphoid} if $Y/E'$ is not ambiguous or need not be remembered). Any holomorphoid of $Y/E'$ which is (tempered) anabelomorphic to $X/E$ (i.e. $Y/E'$ satisfies condition {\bf(4)} above)  will be referred to as an \textit{anabelomorphic holomorphoid of $X/E$}. 
\end{defn}
\brem 
Each anabelomorphic holomorphoid $(Y/E',(E'\into K, K^\flat\isom F), *_{K}:\sM(K)\to \yan_{E'})$ of $X/E$, provides analytic spaces $\yan_{E'}$ and $\yan_K$ a base change morphism $\yan_K\to\yan_{E'}$ (using the given embedding $E'\into K$) which are obtained by Berkovich analytification of the schemes $Y/E'$, $Y\times_{E'}\Spec(K)/K$. The analytification $\yan_K$ must be thought of as the geometric analytic or holomorphic  functions provided by this holomorphoid. This is consistent with the classical theory of analytic functions--where one works with analytic functions over an algebraically closed field.
\erem
	
The following is immediate from the above defintion:

\blem Notably $(X/E,(E\into K, K^\flat\isom F), *_{K}:\sM(K)\to \xan_{E})$ is an anabelomorphic holomorphoid of $X/E$ and  the arithmetic Teichmuller space $\fjxe$ is the space of all anabelomorphic holomorphoids of $X/E$. 
\elem

\para I will often abbreviate  the notation for the holomorphoid $(Y/E',(E'\into K, K^\flat\isom F), *_{K}:\sM(K)\to \yan_{E'})$  as $(Y/E',E'\into K)$.  Note that this short notation $(Y/E',E'\into K)$ is for compatibility with the notational conventions of \cite{joshi-teich}, but let me emphasize that in this paper the tilting data provided by $(Y/E',(E'\into K, K^\flat\isom \cpt), *_{K}:\sM(K)\to \yan_{E'})$ will be very important for computing local arithmetic degrees and whenever needed, I will revert to the long-hand notation described above. 
The properties of $\fjxe$ (resp. $\fjxe_{\cpt}$) are described in \cite[\ssep 1.4]{joshi-teich} and \cite{joshi-untilts}.

\para The case of interest is this: let $E$ be a $p$-adic field, $F=\cpt$
and let $C/E$ be  an elliptic curve. Now  let me indicate how these considerations can be applied to the case of hyperbolic curves considered in \iut. Let $X=C-\{O\}$ be  the canonically punctured elliptic curve which I will simply refer to as the \textit{elliptic cyclops} associated to $C/E$.  Observe that $X/E$ is an hyperbolic curve of topological type $(1,1)$ over the $p$-adic field $E$. 

Later on I will  assume that $C/E$ arises from a number field and hence $X/E$ is  a once punctured elliptic curve i.e. an elliptic cyclops defined over a number field. Such an $X/E$ is of strict Belyi type in the sense of \topics. In that case (i.e. $X/E$ of strict Belyi Type) as was noted in \cite{joshi-teich} (and proved in \cite{mochizuki05}), for any anabelomorphic holomorphoid $(Y/E',E'\into K)$ of $X/E$, one has an isomorphism of $\Z$-schemes $Y\isom X$ i.e. $\fjxe$ is connected in the sense of \cite{joshi-untilts}.

Note that if $(X/E,(E\into K, K^\flat\isom \cpt), *_{K}:\sM(K)\to \xan_{E})$ is a holomorphoid of $X/E$, then the composite morphism $*_{K}:\sM(K)\to\xan_E\to C^{an}_{E}$ provides a geometric base point and hence one obtains a holomorphoid $(C/E,(E\into K, K^\flat\isom \cpt), *_{K}:\sM(K)\to C^{an}_{E})$. Conversely if $*_{K}:\sM(K)\to C^{an}_{E}$ is a basepoint with its image factoring through $\xan_E\to C^{an}_E$ then one obtains from a holomorphoid of $C/E$ a holomorphoid of $X/E$. I will habitually use this to make transition from holomorphoid of $X/E$ to those of $C/E$ and vice versa.

\para Of special interest to us is the situation when $C/E$ has potentially split multiplicative reduction. After a finite extension, this is the context of Tate elliptic curves for which \cite{silverman-arithmetic}, \cite{roquette-book} are  convenient references.

Let $(C/E,(E\into K, K^\flat\isom \cpt), *_{K}:\sM(K)\to C^{an}_{E})$ be an arithmetic holomorphic structure on $C/E$ and that $C/E$ has potentially multiplicative reduction. Since $K/E$ is an algebraically closed extension of $E$, by standard results on Tate elliptic curves (see \cite{silverman-arithmetic}),  there exists some finite extension $E\into E_1\into K$ such that $C/E_1$ has split multiplicative reduction and may be described as a Tate elliptic curve over $E_1$. In particular $C/K$ is a Tate elliptic curve and hence has a Tate parameter which is an element of $K^*$.

Hence one may associate to the holomorphoid $(C/E,(E\into K, K^\flat\isom \cpt), *_{K}:\sM(K)\to C^{an}_{E})$, a Tate parameter which will be denoted by $$q_{(C/E,(E\into K, K^\flat\isom \cpt), *_{K}:\sM(K)\to C^{an}_{E})} \in K^*$$ (or by  $q_{(C/E,(E\into K, K^\flat\isom \cpt))}$ or more compactly $q_{C;K}$ to indicate its dependence on the arithmetic holomorphic structure) and call this the Tate parameter of  the holomorphoid $(C/E,(E\into K, K^\flat\isom \cpt), *_{K}:\sM(K)\to C^{an}_{E})$.  \textit{Note that $C/E$ may not be a Tate curve--so there is some abuse of terminology here!} At any rate the Tate parameter can be viewed as a function (in the sense this term is used in algebraic geometry) on $\fjxe$.

The following is a fundamental observation in understanding how the theory of arithmetic holomorphic structures dictates that local height contributions have dependency on arithmetic holomorphic structures. A similar assertion is tacitly made and used in \iut.

\bthm\label{th:tate-parameter-as-a-function} Suppose that $C/E$ is an elliptic curve over a $p$-adic field $E$ and $X=C-\{O \}/E$ is the canonical elliptic cyclops associated to $C/E$. 
\benumlab
\item If $X/E$ has potentially multiplicative reduction, then every anabelomorphic holomorphoid of $X/E$ i.e. every object of $\fjxe$ has this property. 
\item Hence the Tate parameter (in the above defined sense) provides a function (in the sense of algebraic geometry) on the $\fjxe$ given as follows: $$(Y/E',(E'\into K, K^\flat\isom F), *_{K}:\sM(K)\to \yan_{E'})\mapsto q_{(Y/E',(E'\into K, K^\flat\isom F), *_{K}:\sM(K)\to \yan_{E'})}\in K^*.$$
\item In particular, by the results of \cite{joshi-teich,joshi-untilts}, as arithmetic holomorphic structures may not be comparable, the Tate parameter is a non-constant function on $\fjxe$  and this is true even if one works with $\fjxe_{\cpt}$.
\eenum
\ethm
\bp 
Let $(Y/E',(E'\into K, K^\flat\isom F), *_{K}:\sM(K)\to \yan_{E'})\in\fjxe$ be an anabelomorphic holomorphoid of $X/E$. By construction of $\fjxe$, $Y/E'$ has $$\dim(Y/E')=\dim(X/E)=1$$ i.e. $Y/E'$ is a geometrically connected smooth, quasi-projective curve over $E'$. 

Let the topological type of $Y/E'$ be $(g,r)$. Then by \cite[Lemma 1.3.9]{mochizuki04} one knows that the topological type of $Y/E'$ is amphoric i.e. the topological type of $Y/E'$ is determined by the anabelomorphism class $Y/E'$. As $Y/E'$ and $X/E$ are anabelomorphic, the topological type $(g,r)$ of $Y/E'$ coincides with the topological type of $X/E$ i.e. $(g,r)=(1,1)$. 

Hence there exists some projective  curve $C'/E'$ of genus one over $E'$ and an $E'$-rational point $P$ such that $Y=C'-\{P\}$ i.e. $Y/E'$ is an elliptic cyclops. Thus one sees that if $X/E$ is an elliptic cyclops then every object of $\fjxe$ arises from some elliptic cyclops over some $p$-adic field which itself is anabelomorphic to the $p$-adic $E$.

By \cite[Theorem 4.12(ii)]{mochizuki-topics1}, any anabelomorphism between $Y/E'$ and $X/E$  induces an isomorphism between their geometric tempered fundamental subgroups (working with $\C_p$ and the algebraic closures $\bEp\subset \C_p$ and $\bE\subset \C_p$  of $E'$ and $E$ respectively in $\C_p$  to compute the geometric tempered fundamental subgroups) $$\pitb{Y/E'}\isom \pitb{X/E}$$
and hence by \cite[Corollary 3.11]{mochizuki-semigraphs}  one sees that this anabelomorphism induces an isomorphism of the dual  graphs of stable reductions of $Y/\bEp$ and $X/\bE$. 

Hence one sees that if $X/\bE$ arises from a Tate elliptic curve then so does $Y/\bE'$. In particular  $\yan_K$ is the analytification of a  Tate curve elliptic curve over $K$ with one point removed. In particular if the Tate parameter makes sense for $X/E$ (after passage to some finite extension) then it make sense for every anabelomorphic holomorphoid of $X/E$ i.e. for every object of $\fjxe$ and hence the Tate parameter is a function on $\fjxe$: 
$$(Y/E',(E'\into K, K^\flat\isom F), *_{K}:\sM(K)\to \yan_{E'})\mapsto q_{(Y/E',(E'\into K, K^\flat\isom F), *_{K}:\sM(K)\to \yan_{E'})}\in K^*.$$
This proves {\bf(1)} and  {\bf{(2)}}. The remaining assertion follows from the results regarding arithmetic holomorphic structures established in \cite{joshi-teich,joshi-untilts}.
\ep

\brem 
A central question which now arises (in my theory as well as in \iut) is how to compare values of the Tate parameter function on $\fjxe$ in one common location. This question is also central to \moccor. Mochizuki's way of doing this is discussed in \cite{mochizuki-iut3} (this depends on his approach to arithmetic holomorphic structures which I have discussed in \cite{joshi-untilts}). A different and a more natural solution to this problem was proposed in \cite{joshi-teich} and that strategy will be elaborated here as it lies at the heart of my approach to \moccor\ and it provides Mochizuki's approach as well. This is elaborated here and in \cite{joshi-teich-rosetta}.
\erem

\section{The choice of a theta function}\label{se:theta-function}
\para Let $(X/E,(E\into K, K^\flat\isom F), *_{K}:\sM(K)\to \xan_{E})$ be a holomorphoid of  $X/E$ such that $X/E$ is a hyperbolic curve of type $(1,1)$ arising from a Tate elliptic curve $C/E$ and that the $p$-adic field $E$ satisfies: $$E=E(\sqrt{-1},\sqrt[2\ell]{q},\sqrt[2\ell]{1})\subset K$$  
for some root $\sqrt[2\ell]{q}\in K$ of the Tate parameter of $C/E$ and some roots $\sqrt{-1},\sqrt[2\ell]{1}\in K$ for some  odd prime $\ell\neq p$ which will be specified later. These considerations can also be applied to any $(Y/E',(E'\into K, K^\flat\isom F), *_{K}:\sM(K)\to \yan_{E'})\in\fjxe$ in place of $(X/E,(E\into K, K^\flat\isom F), *_{K}:\sM(K)\to \xan_{E})$ and will be left to the readers.

In \cite{mochizuki-iut3} this hypothesis on $E$ is available through the choice of \textit{Initial $\Theta$-Data} in \iut.  

\textit{To avoid extremely cumbersome notation, in the subsequent paragraphs the data of the holomorphoid $(X/E,(E\into K, K^\flat\isom F), *_{K}:\sM(K)\to \xan_{E})$ will be notationally suppressed, but as should be clear from Theorem~\ref{th:tate-parameter-as-a-function} it is important to remember the dependence on the choice of the holormorphoid.}

\brem To understand the dependence of function theory  on the choice of holomorphoids of $X/E$ (established in \cite{joshi-teich,joshi-untilts}) the following will be useful. 

Let $(X/E,(E\into K, K^\flat\isom F), *_{K}:\sM(K)\to \xan_{E})$ be a holomorphoid of $X/E$.  Then this choice provides the analytic space $C^{an}_K$ which is the analog of a complex structure in the theory of Riemann surfaces.  Suppose that our  $X/E$ arises from  a Tate elliptic curve $C/E$. Then the function theory on $C_K^{an}$ is  described by a theory of $p$-adic $\theta$-functions with values in $K$ (by \cite{roquette-book}).  

Suppose now that one has two distinct holomorphoids of $X/E$ providing perfectoid fields $K_1,K_2$. By \cite[Theorem~3.15.1]{joshi-teich} one knows that if $K_1,K_2$ are two algebraically closed perfectoid fields which are not topologically isomorphic then  the analytic function theories of $C_{K_1}^{an}$ and $C_{K_2}^{an}$ are not isomorphic (as these two analytic spaces are not isomorphic).  Suppose one chooses a  $\theta_M$  as a $\theta$-function over $C^{an}_K$. Then such a function can be evaluated at the elements of $E$ (or finite extensions of $E$). This allows us to talk about values of the $\theta$-functions, even at elements  of $E$ (a finite extension of $\Q_p$), as living in two topologically distinct perfectoid fields. Here I will describe how one can compare these values  based on the idea of \cite{joshi-teich}.
\erem 

\para I will  choose a theta   function $\theta_M$ (described below) on $C/E$. For the theory of $p$-adic theta-functions see \cite{roquette-book}.  A theta-function on $C$ is a function on the  universal cover with some quasi-periodicity properties--equivalently a theta-function is a section of some line bundle.
Let $q=q_E$ be the a Tate parameter for $C/E$. 

Let $(X/E,(E\into K, K^\flat\isom F), *_{K}:\sM(K)\to \xan_{E})$ be a holomorphoid of $X/E$. 
Let $\iota:E\into K$ be the given embedding of valued fields. The chosen theta function $\theta_M$ can be considered as a theta function in $C_K^{an}$. In particular the theta function $\theta_M$ can be viewed as a $\theta$-function over $C_K^{an}$. To avoid complexity of the notation, I will use $\theta_{M;K}$ as a mnemonic for remembering the chosen holomorphoid $(X/E,(E\into K, K^\flat\isom F), *_{K}:\sM(K)\to \xan_{E})$ of $X/E$.

The theta function  I choose here, denoted by  $\tm u$   is the same as the function used in \cite[Proposition 1.4]{mochizuki-theta}  and \iut\ ($M$ in the subscript is an abbreviation of Mochizuki):
\be
\tm u = q_E^{-1/8}\sum_{n\in\Z}q_E^{\frac{1}{2}(n+\frac{1}{2})^2}u^{2n+1}.
\ee
This satisfies the following properties
\benumlab
\item $\tm{u^{-1}}=-\tm u$,
\item $\tm {q_E^\frac{j}{2}u}=(-1)^j q_E^{-\frac{j^2}{2}}u^{-2j}\tm u$.
\eenum
Since one is working with $\sqrt{q_E}$ as opposed to $q_E$, this theta function naturally lives on a double cover $ C'\to C$ of $C$.

From the first formula one sees that $$\tm 1 =0,$$ using this with the quasi-periodicity given by (2) one sees that
$$\tm {q_E^\frac{j}{2}}=0 \text{ for all } j\in \Z,$$
and that $\tm u=0$ if and only if $u=q_E^{j/2}, j\in\Z$.
From this one can also provide a multiplicative description of this theta function along the lines of \cite{roquette-book} or \cite{silverman-advanced}. Secondly the points $q_E^{j/2}$ $(j\in\Z)$ are mapped to the identity element $O\in C'$ under the Tate parametrization of $C'$.

\para There is a line bundle $\sL_\theta$ on $C'$ corresponding to $\tm u$. This line bundle may be trivialized on $Y=C'-\{O\}$ and one sees from the above that  $\tm{u}$ has no zeros on $Y$ and can be used to trivialize $\sL_\theta$ on $Y$. I will colloquially speak of the $\theta$-function $\tm{u}$ as a function on $Y$ (and by abuse of terminology, sometimes even on $C'$). This ties up with \cite[\ssep 1.2]{joshi-teich}.

\section{Choice of theta values}\label{se:theta-values}
\para  Let $\ell\geq 3$ be a prime with $\ell\neq p$. The choice of $\ell$ will be made later on. \textit{Assumptions of \ssep\ref{se:theta-function} will now be in force.}  Let $(X/E,(E\into K, K^\flat\isom F), *_{K}:\sM(K)\to \xan_{E})$ be a holomorphoid of $X/E$. 

Let $\bE\subset K$ be the algebraic closure of $E$ in $K$ (recall that $K$ is an algebraically closed perfectoid field). Let $\zeta_\ell\in K$ be a primitive $\ell^{th}$-root of unity.  For $\alpha=0,\ldots,\ell-1$, let $q^{1/2\ell}\zeta_\ell^\alpha\in K$ be all the $\ell^{th}$-roots of $\sqrt{q}$ in $K$. 

Then the following collection of elements will be referred to as $\theta$-values (for the given holomorphoid data):
$$\frac{1}{\xi_j}=\frac{\tm {{q^\frac{j}{2\ell}\zeta_\ell}}}{\tm {\zeta_\ell}}=(-1)^j {q^{-\frac{j^2}{2\ell}}}\zeta_\ell^{-2j} \in K.$$
Note that one has
$$\xi_j\subset\O_{\bE}\subset \O_K.$$
The theta values of interest to us are the $\theta$ values $\xi_j$ for $j=1,\ldots,\ells$. There is an evident action of $G_E$ on the set of values considered above. 

By the description of the torsion points of a Tate curve in \cite{serre-abelian}, the images of $q^{1/2\ell}\zeta_\ell^\alpha$ under the Tate parameterization of $C'$ provide a subgroup of order $\ell$ in $C'[\ell]$ on the double cover of $C'\to C$. 

\para Note that I will write $$\theta_{M;K}(q^{1/2\ell}\zeta^j_{2\ell})$$ and 
$$\frac{1}{\xi_{j;K}}=\frac{\theta_{M;K} ({{q^\frac{j}{2\ell}\zeta_\ell}})}{\theta_{M;K} ({\zeta_\ell})}$$  to indicate the dependence on the holomorphoid $(X/E,(E\into K, K^\flat\isom F), *_{K}:\sM(K)\to \xan_{E})$ of $X/E$. I hope that my notational simplification $\xi_{j;K}$ will not be too confusing.

\section{Mochizuki's ansatz ``$q=q^{j^2}$'' and the $\Theta_{gau}$-links}\label{se:ansatz}
\para Let $\ell\geq 3$ be a prime number with $\ell\neq p$ and let $\ells=\frac{\ell-1}{2}$. Let $q$ be the Tate parameter  of a Tate elliptic curve over a $p$-adic field. According to \cite[Page 4]{mochizuki-gaussian}, Mochizuki's theory of $\theta$-links  is roughly based on the idea that at some level in his theory \iut\ one can declare that   (I am putting this in quotes intentionally) $$\text{``}q=q^{j^2}\text{''} \text{ and } \abs{q}_E<1.$$ 
In \iut\ Mochizuki, in fact, asserts that in his theory one can set
$$\text{``}\left(q^{1^2},q^{2^2},q^{3^2},\cdots, q^{\ells^2}\right)=q\text{''}.$$
This  assertion is the central part of the theory of $\Theta_{gau}$-links in Mochizuki's theory  and is central to the construction of $\ttheta_{Mochizuki}$.  My construction is detailed in \ssep\ref{se:construction-ttheta}. 

It is obviously enough to  be able to set
\be
\label{eq:ansatz-for-p}\text{``}p=\left(p^{1^2},p^{2^2},p^{3^2},\cdots, p^{\ells^2}\right)\text{''}.
\ee

\para \textit{The important discovery of the present paper is that this can be done canonically and explicitly in the Arithmetic Teichmuller Theory of \cite{joshi-teich} as I will now demonstrate.} This paper provides an intrinsic geometric approach to the existence of Mochizuki's theory of the $\Theta_{gau}$-Links \cite{mochizuki-iut3}.

Let $F$ be an algebraically closed, perfectoid field of characteristic $p>0$ with its given valuation $\abs{-}_F$. Typical example of interest is $F=\cpt$, but much of the theory detailed here works in general and is needed for working with $\fjxe$ (as opposed to $\fjxe_{\cpt}$).  The following proposition will be used throughout.

\bpro\label{pr:primitive-elements1}
Fix $0\neq a\in\fm_F$. Then for $j=1,\ldots,\ells$ 
\benumlab
\item one has primitive elements of degree one   given by $$[a^{j^2}]-p\in W(\O_F).$$ 
\item Each of these elements generates a principal prime ideal  
$$\wp_j=([a^{j^2}]-p)\subset W(\O_F).$$ 
\item Hence each  $\wp_j$ defines  a closed classical points $y_{j}\in \syfQp$.
\item For each $y_j$ one  may  consider its image  under the canonical morphism $\syQp\to\sxQp$ and view $y_j$ as providing a closed classical point of $\sxqp$. 
\eenum
\epro
\bp 
The first assertion is immediate from \cite[D\'efinition 2.2.1]{fargues-fontaine}. The second assertion is immediate from \cite[Section 2.2]{fargues-fontaine}. By \cite[Lemma 2.2.14]{fargues-fontaine} for each $1\leq j\leq \ells$, the ideal $$\wp_j=([a^{j^2}]-p)\subset W(\O_F)$$ is a principal prime ideal of $W(\O_F)$ generated by a primitive element of degree one and its extension to $W(\O_F)\subset B$ provides a closed maximal ideal of $B$ and hence a closed classical point of $\syfQp$. The remaining assertion is clear.
\ep

\blem\label{pr:primitive-elements2} 
Let $\wp\subset W(\O_F)$ be a principal (prime) ideal generated by some primitive element of degree one. Then there is a tuple of principal prime ideals $(\wp_1,\wp_2,\ldots,\wp_\ells)$ of the  sort  constructed in Proposition~\ref{pr:primitive-elements1} with $\wp=\wp_1$. 
\elem 
\bp 
Let $\wp=(\alpha)\subset W(\O_F)$ be a principal ideal generated by a primitive element $\alpha$ of degree one. Then by \cite[Corollaire 2.2.9]{fargues-fontaine} one can find a unit $u\in W(\O_F)^*$ and an $a\in\fm_F-\{0\}$ such that $$\alpha\cdot u=[a]-p$$
and by definition $[a]-p$ is a primitive element of degree one in $W(\O_F)$.
Hence $$\wp=(\alpha)=([a]-p)$$ and hence  one has the $\ells$-tuple of primitive elements of degree one generating prime ideals $\wp_j=([a^{j^2}]-p)$ for $j=1,\ldots,\ells$ with $\wp=\wp_1$ as asserted.
\ep
 
\begin{defn}\label{def:primitiv-ansatz}
	Let $F$ be an algebraically closed perfectoid field of characteristic $p>0$. Then \textit{Mochizuki's Primitive Ansatz},  is the subset
	$$\tsigf=\left\{ ([a]-p,[a^{2^2}]-p,[a^{3^2}]-p, \cdots ,[a^{\ells^2}]-p): a\in\fm_F-\{0\}\right\}$$
	 of the set of $\ells$-tuples of primitive elements of degree one of $W(\O_F)$.
\end{defn}

\para
I will usually view an element of $\tsigf$ as providing a tuple of points of $\syfQp$ i.e. I will usually conflate $[a^{j^2}]-p$ with the corresponding closed classical point $y_j$ in $\syfQp$. Hence each element of $\tsigf$ provides an $\ells$-tuple of closed classical point of $\syfQp$ i.e. $\tsigf$ may be identified with its image in  $$\tsigf\subset\abs{\syfQp}^\ells=\overbrace{\abs{\syfQp}\times \abs{\syfQp} \times \cdots \times \abs{\syfQp}}^{\ells\text{ factors}}.$$  
\textcolor{red}{Note that an arbitrary tuple in $\syfQp^{\ells}$ need not belong to $\tsigf$.} So the set $\tsigf$ should not be confused with $\syfQp^\ells$.

\bcor
$\tsigf$ viewed as a subset of $\abs{\syfQp}^\ells$  is metrisable.
\ecor
\bp 
This is immediate from the fact that the set of closed classical points $\abs{\syfQp}$ is a metric space by \cite[Proposition 2.3.2]{fargues-fontaine} and hence $\abs{\syfQp}^\ells$ is metrisable and so $\tsigf$ may be equipped with a metric induced from $\abs{\syfQp}^\ells$.
\ep 

\para Each  $([a]-p,[a^{2^2}]-p,[a^{3^2}]-p, \cdots ,[a^{\ells^2}]-p)\in\tsigf$ defines an $\ells$-tuple of closed classical points $(y_1,\ldots,y_\ells)\in \syfQp^\ells$ and hence also provides us with tuple of untilts of $F$
$$(K_j,K_j^\flat\isom F)_{1\leq j\leq\ells}=(K_{y_j}, K_{y_j}^\flat\isom F)_{1\leq j\leq\ells}$$ where $K_j$ are the residue fields $$(K_j)_{1\leq j\leq\ells}=(K_{y_j})_{1\leq j\leq\ells}$$ each of which is an algebraically closed, perfectoid and of characteristic zero.

\para Notably, let $t\in \cpt$ be such that $([t]-p)\subset W(\O_{\cpt})$ gives the canonical point of $\sxQp$ (\cite[Chap 10]{fargues-fontaine}) with residue field $\C_p$ as a $\gal(\bQ_p/\Q_p)$-module.  Then the tuple
\be\label{eq-canonical-tuple}
([t]-p,[t^{2^2}]-p,[t^{3^2}]-p, \cdots ,[t^{\ells^2}]-p)\in\tsigcpt
\ee
in which the first term  corresponds to the canonical point of $\syQp$. 

\para Now let $$\vphi:\syfQp\to\syfQp$$ be the Frobenius morphism of $\syfQp$ (\cite{fargues-fontaine}). For $F=\cpt$,  the fiber of $\syQp\to\sxQp$ over the canonical point of $\sxQp$ can be identified with the set of prime ideals generated by primitive irreducible elements of degree one:
$$
\left\{([\vphi^n(t)]-p)\subset W(\O_{\cpt}):n\in\Z \right\}.
$$
For later use let me record the following assertion:
\bpro\label{pr:frob-stable-cannonical}
Let $F$ be an algebraically closed, perfectoid field of characteristic $p>0$. Let $0\neq a\in \fm_F$ and let $([a]-p)$ be a principal prime ideal of degree one. Then 
\benumlab 
\item For each $n\in\Z$, one has
$$([\vphi^n(a)]-p,[\vphi^n(a^{2^2})]-p,[\vphi^n(a^{3^2})]-p, \cdots ,[\vphi^n(a^{\ells^2})]-p)\in\tsigcpt.$$ 
\item Hence $\tsigf$ is stable under iterates of Frobenius $\vphi:\syfQp\to \syfQp$. 
\item 
In particular $\tsigcpt$ contains the Frobenius orbit of the canonical tuple \eqref{eq-canonical-tuple}.
\eenum
\epro
\bp 
Clearly ${\bf(1)}\implies{\bf(2)}\implies{\bf(3)}$. So it is enough to establish {\bf(1)}. Observe that for $j=1,\ldots,\ells^2$ one has
$$[\vphi^n(t^{j^2})]-p=[\vphi^n(t)^{j^2}]-p$$  and so if I write $b=\vphi^n(a)$ then for $j=1,\ldots,\ells$ one has
$$[\vphi^n(a^{j^2})]-p=[b^{j^2}]-p$$ and hence the tuple
$$([b]-p,[b^{2^2}]-p,[b^{3^2}]-p, \cdots ,[b^{\ells^2}]-p)$$
is contained in $\tsigf$ by the definition of $\tsigf$. This proves {\bf(1)}.
\ep 

\brem

One should think of the assignment $$([a]-p)\mapsto ([a^{1^2}]-p), ([a^{2^2}]-p), \ldots, ([a^{\ells^2}]-p)$$ 
as divisorial correspondence $\syQp\dashrightarrow \syQp$ (resp. $\sxQp\dashrightarrow\sxQp$) in the sense of a Hecke correspondence between modular curves. 
\erem

\para Now suppose that $F=\cpt$ equipped with its natural Galois action of $G_{\Q_p}=\gal(\bQ_p/\Q_p)$ which is studied in \cite[Chapitre 10]{fargues-fontaine}. There is a natural continuous action of $G_{\Q_p}$ on $\syQp$ via the natural $G_{\Q_p}$-action on the set of primitive 
elements of degree one of $W(\O_{\cpt})$. On a tuple $(y_1,\ldots,y_{\ells})\in\tsigcpt$ this action is given by
$$y_j\mapsto \sigma(y_j)=\sigma([a^{j^2}]-p)=[\sigma(a)^{j^2}]-p$$
so that $(\sigma(y_1),\ldots,\sigma(y_{\ells}))\in\tsigcpt$. 
Note that one has the mapping 
$$\tsigcpt\to \syQp \to \sxQp$$
given by $$(y_1,\ldots,y_{\ells})\mapsto y_1\to \text{ image}(y_1)\in \sxQp.$$
Hence one has proved the following proposition:
\bpro\label{pr:prim-ansatz-galois-stable}
Mochizuki's Primitive Ansatz $\tsigcpt$ is Galois stable (and Frobenius stable by Proposition~\ref{pr:frob-stable-cannonical}), and notably any point in the Frobenius orbit of the canonical tuple of \eqref{eq-canonical-tuple} has the property that the image of  its first coordinate in $\sxQp$ is the canonical point of $\sxQp$ and hence the image in $\sxQp$ is fixed  under the action $G_{\Q_p}$.
\epro

\para\label{pa:norms-vs-element-comp} Suppose that $E$ is a $p$-adic field and   $K_1,\ldots,K_\ells$ are algebraically closed perfectoid, each equipped with an isometric embedding $\iota_j:E\into K$. Then as each $K_j$ is algebraically closed so  each field $K_j$ contains a copy of the algebraic closure of $E$,  $\bE_j\subset K_j$. As was observed in \cite{joshi-teich}  there may be no topological isomorphism $K_1\isom K_2$ (say) which takes $\iota_1:\bE_1\into K_2$ isomorphically to $\iota_2:\bE_2\into K_2$.   Notably one cannot compare elements $\iota_j(z)\in K_j$ with each other for any $z\in E$ (or  $z\in \bE$). However  because of \ssep\ref{pa:value-group-comp-loc} one may be able to compare the absolute values $\abs{\iota_1(z)}_{K_1}$ and $\abs{\iota_2(z)}_{K_2}$. I will simply write $\abs{z}_{K_1}$ instead of $\abs{\iota_1(z)}_{K_1}$. Hopefully this will not cause any confusion. Thus one must take care when working with elements of our field $E$ and its algebraic closures in the fields $K_j$.

\para Theorem~\ref{pr:lift-vals} demonstrates a crucial property of points  of Mochizuki's Primitive Ansatz \ref{pa:Mochizuki-Ansatz}. This property shows that any point $(y_1,\ldots,y_\ells)\in\tsigf$ of Mochizuki's Primitive Ansatz provides a corresponding tuple of perfectoid fields $(K=K_1,K_2,\ldots,K_\ells)$ whose valuations are simultaneously scaled by suitable factors--\emph{this simultaneous scaling of valuations is the primary reason why this ansatz  is  important in \iut.}
\bthm\label{pr:lift-vals} 
Let $(y_1,\ldots,y_\ells)\in\tsigf$ be viewed as a tuple of closed classical points of $\syfQp$. Then one has (remembering \ssep\ref{pa:norms-vs-element-comp})
\benumlab
\item $$v_{K_j}(p)=j^2v_{K_1}(p)\text{ for } j=1,\ldots,\ells.$$
\item Notably if $K_1=\C_p$ then $v_{K_j}(p)=j^2$.
\item Hence for  $z\in \bE\subset K_j$, for $j=1,\ldots,\ells$ one has
$$\abs{z}_{K_j}=\abs{z}_{K_1}^{j^2}.$$
\item However the valued fields $K_1,\ldots,K_\ells$ need not be all topologically isomorphic.
\eenum
\ethm

\bp 
\textcolor{red}{Note the since one can view the values $\abs{-}_{K_j}$ in the value group of $F$  by \ssep\ref{pa:value-group-comp-loc} one can make such assertions about $\abs{-}_{K_j}$.}

The valuation of $K_j$ can be computed from \cite[Proposition 2.2.17]{fargues-fontaine} as $$v_{K_j}(p)=v_F(a^{j^2})=j^2v_F(a)=j^2v_{K_1}(p).$$
Hence the first assertion is proved and the second follows from this. Note that by loc. cit. the restriction of $\abs{-}_{K_j}$ to $E\subset K_j$ is given by the above formula and hence the penultimate assertion is also established. The last claim is an immediate  consequence of \cite{kedlaya18}.
\ep

\para\label{pa:Mochizuki-Ansatz} Let $(y_1,\ldots,y_\ells)\in\tsigf$ be a point of Mochizuki's Primitive Ansatz and let me write $K$ for the residue field $K_{y_1}$ of $y_1$. So one has  $K\supset \Q_p$ and $K^\flat=F$. By \cite[Corollaire 2.2.9]{fargues-fontaine} one can choose a primitive element of degree one of $([p^\flat]-p)\subset W(\O_F)$ which generates the prime ideal $\ker(\eta_K:W(\O_F)\onto \O_K)$ defining the point $y_1$ (and hence $K$) i.e. $([p^\flat]-p)=\ker(\eta_K:W(\O_F)\onto \O_K)$. In this notation the point of  $(y_1,\ldots,y_\ells)\in\tsigf$ is given by equations 
\begin{align*}
[p^{\flat}]-p&=0,\\
[(p^{\flat})^{2^2}]-p&=0\\
\vdots &\quad\ \vdots\\
[(p^{\flat})^{\ells^2}]-p&=0.
\end{align*}
And hence one has
$$p=[p^{\flat}]=[(p^{\flat})^{2^2}]=
\cdots=[(p^{\flat})^{\ells^2}].
$$
This proves the following lemma:
\blem 
Hence every point of $\tsigf$ provides a version of Mochizuki's ansatz \eqref{eq:ansatz}, \eqref{eq:ansatz2}  in the form stated in \eqref{eq:ansatz-for-p}. Explicitly 
\be p\mapsto ([p^{\flat}],[(p^{\flat})^{2^2}],
\cdots,[(p^{\flat})^{\ells^2}]).
\ee
Moreover for each point of $\tsigf$ one has the valuation scaling property given by Theorem~\ref{pr:lift-vals}. 
\elem

\brem 
Let me remark that I use the projection to the first coordinate $$\tsigcpt \subset \syQp^\ells \to\syQp$$ for valuation computations. There are also more complicated (and highly non-algebraic) maps from $\tsigf\to \syQp$ possible for example $$([a]-p,[a^{2^2}]-p,[a^{3^2}]-p, \cdots ,[a^{\ells^2}]-p)\mapsto [a^{\sum_{j=1}^\ells j^2}]-p.$$
\erem

\brem 
The set $\tsigcpt$ satisfying these properties should be viewed as providing arithmetic Kodaira Spencer classes. From \cite[proof of Theorem 3.15.1]{joshi-teich}  one sees that these are non-trivial because as one moves over $\tsigcpt$ using the action of the above groups,  the analytic spaces $\xan_{K_1},\ldots,\xan_{K_\ells}$ need not be topologically isomorphic in general.
\erem

\para\label{pa:Ansatz-in-teich-space}
Let me now demonstrate that the existence of Mochizuki's Primitive Ansatz $\tsigcpt$ (Definition~\ref{def:primitiv-ansatz}) has consequence for arithmetic Teichmuller spaces $\fjxe$. Namely Mochizuki's Primitive Ansatz leads to the construction of $\Theta_{gau}$-Links of \cite{mochizuki-iut3}.  This rests on the following construction.

\begin{defn}\label{def:ansatz}
Let $E$ be a $p$-adic field, let $X/E$ be a geometrically connected, smooth quasi-projective variety over $E$ and let $\fjxe$ be the arithmetic Teichmuller space constructed in \cite{joshi-teich,joshi-untilts}. 
\textit{Mochizuki's Ansatz} is the full subcategory  $\sigfjxecpt$ of the product $\fjxe_{\cpt}^\ells$  whose  set of objects consists of $\ells$-tuples of arithmetic holomorphic structures (see \cite[Section 4]{joshi-untilts}) of the form
\be\label{eq:theta-link-tuple-def}(X/E,(E\into K_{y_j}, K_{y_j}^\flat\isom \cpt), *_{K_{y_j}}:\sM(K_{y_j})\to \xan_{E}))_{j=1,\ldots,\ells}\in\fjxe_{\cpt}^\ells,\ee
where $(y_1,\ldots,y_\ells)\in\tsigf$	 i.e. $\sigfjxecpt$ is the set of $\ells$-tuples of anabelomorphic holomorphoids of $X/E$ lying over some point $(y_1,\ldots,y_\ells)\in\tsigf$ of Mochizuki's Primitive Ansatz $\tsigcpt$.
\end{defn}

\brem 
One can also work with a wider definition of Mochizuki's Ansatz. This will be the full subcategory whose object set is the $\sigfjxe\subset\fjxe^\ells$ consisting of the objects of the form 
$$(X/E,(E\into K_{y_j}, K_{y_j}^\flat\isom F), *_{K_{y_j}}:\sM(K_{y_j})\to \xan_{E}))_{j=1,\ldots,\ells}\in\fjxe_{F}^\ells$$
lying over a point $(y_1,\ldots,y_\ells)\in\tsigf$ for some algebraically closed perfectoid field $F$ of characteristic $p>0$. However the case $F=\cpt$ is adequate for the first readings.
\erem

\para The following remarks will help the readers in understanding  $\sigfjxe$ and its geometric significance.
\label{re:ansatz-as-correspondence}\
\subpara From my point of view, a point of Mochizuki's Ansatz given by \eqref{eq:theta-link-tuple-def} is a divisorial correspondence on Arithmetic Teichmuller Space $\fjxe_{\cpt}$ quite similar to the classical Hecke correspondences on modular curves. Especially  one should write \eqref{eq:theta-link-tuple-def}  in the divisorial notation as formal sum of ``points'' of the arithmetic Teichmuller space $\fjxe_{\cpt}$.
$$\sum_{j=1}^\ells[(X/E,(E\into K_{y_j}, K_{y_j}^\flat\isom \cpt), *_{K_{y_j}}:\sM(K_{y_j})\to \xan_{E}))].$$
\subpara One wants to apply this divisorial correspondence to functions on $\fjxe$ (such as the one provided by Theorem~\ref{th:tate-parameter-as-a-function}) in a manner similar to averaging of modular functions over a Hecke correspondence. \textit{This is the core idea behind (and especially my approach to) \moccor.}
\subpara Notably  one wants to replace the formal sum in \eqref{eq:theta-link-tuple-def} by   the sum over the values of the ``Tate parameter function'' given by Theorem~\ref{th:tate-parameter-as-a-function}:
$$
\sum_{j=1}^\ells q_{(X/E,(E\into K_{y_j}, K_{y_j}^\flat\isom \cpt), *_{K_{y_j}}:\sM(K_{y_j})\to \xan_{E}))}.
$$
since the Tate parameter is multiplicative it makes more sense to write this multiplicatively:
$$
\prod_{j=1}^\ells q_{(X/E,(E\into K_{y_j}, K_{y_j}^\flat\isom \cpt), *_{K_{y_j}}:\sM(K_{y_j})\to \xan_{E}))}.
$$

\subpara \textit{However} as was observed in \cite{joshi-teich}, there is no common topological field in which one can view the individual terms of the above sum or product. \textit{Hence neither this sum nor the product  makes sense.} This fundamental problem is resolved in \ssep~\ref{se:lifting-values-to-B}.

\para Objects of $\sigfjxecpt$ are pointed objects (i.e. equipped with a distinguished point). Often it is convenient to forget the base-point information i.e. there is a ``forget-all-the-base-points'' functor which maps
$$(X/E,(E\into K_{y_j}, K_{y_j}^\flat\isom \cpt), *_{K_{y_j}}:\sM(K_{y_j})\to \xan_{E}))_{j=1,\ldots,\ells}\in\fjxe_{\cpt}^\ells$$
to
$$(X/E,(E\into K_{y_j}, K_{y_j}^\flat\isom \cpt))_{j=1,\ldots,\ells}$$
(such tuples are a category in the obvious way).

The following is now clear from Proposition~\ref{pr:frob-stable-cannonical} and Proposition~\ref{pr:prim-ansatz-galois-stable} and the above discussion (also see \cite{joshi-teich}):
\bpro 
There is a natural action of $G_{\Q_p}$, and Frobenius $\vphi$ of $\syQp$ on the tuples
$$(X/E,(E\into K_{y_j}, K_{y_j}^\flat\isom \cpt))_{j=1,\ldots,\ells}$$
given as follows. For any $\sigma\in G_{\Q_p}$ or $\sigma\in\vphi^\Z$ the object
$$(X/E,(E\into K_{y_j}, K_{y_j}^\flat\isom \cpt))_{j=1,\ldots,\ells}$$
is mapped to
$$(X/E,(E\into K_{\sigma({y_j})}, K_{\sigma({y_j})}^\flat\isom \cpt))_{j=1,\ldots,\ells}.$$
\epro

\section{Lifting theta values to $B$--Existence of Teichmuller lifts}\label{se:lifting-values-to-B}
\para I will now demonstrate how the problem discussed in Remark~\ref{re:ansatz-as-correspondence}{\bf(4)} can canonically resolved, not just for the Tate parameter (Theorem~\ref{th:tate-parameter-as-a-function}) but also for $\theta$-values (\ssep\ref{se:theta-values}). This section is an elaboration of the central idea from \cite{joshi-teich} of lifting elements to the Fargues-Fontaine ring $B$.

\para To understand the lifts of these values (\ssep\ref{se:theta-values}) to $B$, I will use the following proposition. Assumptions on $X/E$ made in \ssep\ref{se:tate-setup} will now be in force notably $X/E$ is a canonical elliptic cyclops arising from a Tate elliptic curve $C/E$.  The formalism of theta-values discussed in \ssep\ref{se:theta-function} and \ssep\ref{se:theta-values} can be applied to any  anabelomorphic holomorphoid $(Y/E',(E'\into K, K^\flat\isom F), *_{K}:\sM(K)\to \yan_{E'}))$ of $X/E$ by the proof of Theorem~\ref{th:tate-parameter-as-a-function} $\yan_K$ arises from a Tate elliptic curve over $K$. For simplicity of notation I will work with the holomorphoids of $X/E$ of the form  $(X/E,(E\into K, K^\flat\isom F), *_{K}:\sM(K)\to \xan_{E}))$. 

\para Let $(X/E,(E\into K, K^\flat\isom F), *_{K}:\sM(K)\to \xan_{E}))$ by any holomorphoid of $X/E$.   For a perfectoid field $K$ with $K^\flat=F$ write $\eta_K:W(\O_F)\to \O_{K}$ for the canonical surjection \cite{fontaine94a} or \cite{fargues-fontaine}. This surjection is traditionally denoted by $\theta$ so the change of notation is naturally forced upon us.

\para For $\rho\in(0,1)\subset \R$ let $\abs{-}_\rho$ be the multiplicative norms defining the Fr\'echet structure of $B$. Let 
\be\label{eq:tate-module-of-lubin-tate} T_K\subset \bpip\subset B\ee be the Tate module of $T_K\subset \sG(\O_K)\isom \bpip$ \cite[Chapitre 4]{fargues-fontaine}
where $\sG/\Z_p$ is a Lubin-Tate formal group over $\Z_p$ with formal logarithm $\sum_{n=0}^\infty\frac{T^{p^n}}{p^n}$. One can also take $\sG=\hgm/\Z_p$.

For a closed point $y$ of degree one of $\syfQp$ with $K=K_y$, I will write $T_y=T_{K_y}$ for simplicity. \emph{Note that $T_y$ is a free $\Z_p$-module of rank one.}

\bpro\label{pr:teichmuller-lift-A}
Let $\xi$ be one of the $\theta$-values (see \ssep\ref{se:theta-values}). Let $y$ be a closed classical point of degree one of $\syfQp$ with residue field $K_y$. Then 
\benumlab 
\item there exists $x\in\fm_F\subset \O_F$ such that 
$\eta_{K_y}([x])=\xi$.
\item $\abs{[x]}_\rho=\abs{x}_F=\abs{\xi}_{K_y}$,
\item and for any $\tau\in T_y\subset \bpip$ one has 
$$\eta_{K_y}(\tau+[x])=\xi.$$
\item So the values $\tau+[x]$ provide a lift of $\theta$-values $\xi$ and notably one has
$$\abs{[x]}_\rho\leq \sup\left\{\abs{\tau+[x]}_\rho:\tau\in T_y\right\}.$$
\eenum
\epro
\bp 
The first assertion has nothing to do with specifics of theta values $\xi$ and is a consequence of \cite[Corollary 2.2.8]{fargues-fontaine} which asserts that there exists a Teichmuller lift of every element of  $ \O_{K}$ under the canonical surjection $\eta_K:W(\O_F)\to\O_K$. The second assertion  just uses the definition of $\abs{-}_{K_y}$ given by \cite[Proposition 2.2.17]{fargues-fontaine} and (1). The third property is immediate from the fact that $\ker(\eta_{K_y})\supset T_y$. In fact, one has the exact sequence of Banach spaces (\cite[Propostion 4.5.14]{fargues-fontaine})
$$\xymatrix{0\ar[r] & T_y\tensor \Q_p\ar[r] &  \bpip \ar[r]^{\eta_{K_y}}& K_y\ar[r] & 0.}$$
The fourth property is self-evident.
\ep

\brem 
Let $[x]$ be a  lift of $\xi$ to $B$ with $\eta_{K_y}([x])=\xi$. Then one is interested in bounds for $\abs{[x]}_\rho=\abs{\xi}_{K_y}$. The idea is to instead consider bounds for elements of $T_y+[x]$.

Note that union of all such lifts $T_y+[x]$ of $\xi$ is not closed under Frobenius. So one needs to enlarge the locus to be Frobenius stable and then consider its intersection with $\bpip$.
\erem

\bpro\label{pr:teichmuller-lift-B}
Let $(X/E,(E\into K_{y_j}, K_{y_j}^\flat\isom \cpt), *_{K_{y_j}}:\sM(K_{y_j})\to \xan_{E}))_{j=1,\ldots,\ells}\in\sigfjxecpt$ be a point of Mochizuki's Ansatz lying over a point of Mochizuki's Primitive Ansatz  $(y_1,\cdots,y_\ells)\in\tsigcpt$. Let $[x_j]\in W(\O_{\cpt})\subset B$ be a Teichmuller lift of $\xi_1$  under
$\eta_{K_j}:W(\O_{\cpt})\to \O_{K_j}$  for $j=1,\ldots,\ells$ given by Proposition~\ref{pr:teichmuller-lift-A}. Then one has
$$\abs{[x_j]}_\rho=\abs{[x_1]}^{j^2}_\rho.$$
\epro
\bp 
This is clear from Proposition~\ref{pr:lift-vals} and the fact that $[x_j]$ is a lift of $\xi_1$ under $\eta_{K_j}$, and hence $$\abs{[x_j]}_\rho=\abs{\xi_1}_{K_j}=\abs{\xi_1}_{K_1}^{j^2}=\abs{[x_1]}_\rho^{j^2}.$$
\ep

\section{Construction of $\ttheta\subset B^\ells$}\label{se:construction-ttheta}
\para\label{pa:existence-theta-pilot-objs} Now one can start the construction of a subset of interest. Let me begin with the following consequence of the results of the previous section. This theorem is the key to the construction of $\ttheta$ detailed below. 

\bthm\label{thm:theta-pilot-object-appears}
Let $E$ and $X/E$ be as in \ssep\ref{se:tate-setup} and \ssep\ref{se:theta-values}. Let $(X/E,(E\into K_{y_j}, K_{y_j}^\flat\isom \cpt), *_{K_{y_j}}:\sM(K_{y_j})\to \xan_{E}))_{j=1,\ldots,\ells}\in\sigfjxecpt$ be a point of Mochizuki's Ansatz lying over  a point $(y_1,\cdots,y_\ells)\in\tsigcpt$ of Mochizuki's Primitive Ansatz and $(\iota_1:E\into K_1,\ldots,\iota_\ells:E\into K_{\ells})$ be the tuple of (perfectoid) residue fields of $(y_1,\cdots,y_\ells)$ with the given isometric embeddings $\iota_j:E\into K_j$. Let $\xi=\xi_1\in E$ be the theta value for $(X/E,\xan/K_1)$ constructed in \ssep\ref{se:theta-values}. Let $$q_j=q_{(X/E,(E\into K_{y_j}, K_{y_j}^\flat\isom \cpt), *_{K_{y_j}}:\sM(K_{y_j})\to \xan_{E})}\in K_j^*$$ be the Tate parameter of the holomorphoid ${(X/E,(E\into K_{y_j}, K_{y_j}^\flat\isom \cpt), *_{K_{y_j}}:\sM(K_{y_j})\to \xan_{E})}$ of $X/E$. 
Let $[x_j]$ be a Teichmuller lift  to $B$ of $q_j$ under $\eta_{K_j}$ defined in Proposition~\ref{pr:teichmuller-lift-B}. Then
\benumlab
\item One has $$q_j=\iota_j(\xi_1)\text{ for } j=1,\ldots,\ells$$ and hence $$\abs{q_j}_{K_j}=\abs{\xi_1}_{K_j}.$$
\item In particular each $[x_j]$ (for $j=1,\ldots,\ells$) is a lift of the theta-value $\xi_1$ to $B$.
\item
\be\label{eq:log-theta-pilot-prod-ver}
\prod_{j=1}^{\ells}\abs{[x_j]}_\rho=\prod_{j=1}^{\ells}\abs{\xi}_{K_1}^{j^2}, 
\ee
\item 
or written additively this is
\be \label{eq:log-theta-pilot}
\sum_{j=1}^{\ells}\log\abs{[x_j]}_\rho=\sum_{j=1}^{\ells}j^2\cdot\log\abs{\xi}_{K_1}.
\ee
\item 
Notably in the left hand side of \eqref{eq:log-theta-pilot},  one can replace the Teichmuller lifts  $\abs{[x_j]}_\rho$ by $\abs{[x_j]+\tau_{j}}_\rho$ where $\tau_j\in T_{y_j}$ and work with the supremum of all such values:
\be 
\sup_{\tau_1,\ldots,\tau_{\ells}}\left\{\sum_{j=1}^{\ells}\log\abs{[x_j]+\tau_j}_\rho\right\} \geq \sum_{j=1}^{\ells}\log\abs{[x_j]}_\rho=\sum_{j=1}^{\ells}j^2\cdot\log\abs{\xi}_{K_1}.
\ee
where the supremum is taken over $\tau_j\in T_{y_j}$ for $j=1,\ldots,\ells$.
\item Notably for $K_1=\C_p$ one has
\be\label{eq:eq:log-theta-pilot-variation} 
\sup_{\tau_1,\ldots,\tau_{\ells}}\left\{\sum_{j=1}^{\ells}\log\abs{[x_j]+\tau_j}_\rho\right\} \geq \sum_{j=1}^{\ells}j^2\cdot\log\abs{\xi}_{\C_p}.
\ee
\item  One has, from the definition of $\xi_1,\ldots,\xi_{\ells}$ that
$$\abs{\xi_j}_{\C_p}=\abs{\xi_1}^{j^2}_{\C_p}$$ and hence one can write 
$$\sum_{j=1}^{\ells}\log\abs{\xi_j}_{\C_p}=\sum_{j=1}^{\ells}\log\abs{\xi_1}^{j^2}_{\C_p}=\sum_{j=1}^{\ells}{j^2}\cdot \log\abs{\xi_1}_{\C_p}.$$
\eenum
\ethm

\bp The first assertion is clear from the properties of the theta functions and theta-values established in \ssep\ref{se:theta-function} and \ssep\ref{se:theta-values}. Moreover the theta function  $\theta_{M;K_j}$ on $\xan_{K_j}$ is the pull-back of the $\theta_{M;E}$ via the morphism of analytic spaces $\xan_{K_j}\to \xan_E$ given by base extension $\iota_j:E\into K_j$. The second assertion follows from  Proposition~\ref{pr:teichmuller-lift-B}. Evidently {\bf(2)}$\implies${\bf(3)}. The fourth assertion follows from Proposition~\ref{pr:teichmuller-lift-A} and obviously {\bf(4)}$\implies${\bf(5)}. This proves the theorem.
\ep

\para\label{re:theta-pilot-obj}\ 
Several important remarks  are in order to put Theorem~\ref{thm:theta-pilot-object-appears} in the perspective and the context of  \moccor.  Notably I want to clarify the relationship with Mochizuki's strategy and the one adopted here. 
\nwsss
\subpara First of all, let me emphasize that Mochizuki was the first to recognize that the value  \eqref{eq:log-theta-pilot}
is a function of arithmetic holomorphic structures, indeterminacies,$\flog$-links and theta-links. My discussion of Arithmetic holomorphic structrues is in \cite{joshi-untilts}. From my point of view indetermincies arise by forgetting arithmetic holomorphic structures (in my sense)--my discussion of this view of Indeterminacies is in \cite{joshi-teich-summary-comments}; my construction of theta-links is \ssep\ref{se:ansatz} and of $\flog$-links is in \ssep\ref{se:log-links}.

This point lies at the core of Mochizuki's construction of $\ttheta_{Mochizuki}$ and hence of \moccor. My Theorem~\ref{thm:theta-pilot-object-appears} provides an independent and conceptually cleaner proof of Mochizuki's assertion based on my approach to arithmetic holomorphic structures detailed in \cite{joshi-teich,joshi-untilts}.
\subpara A tuple 
\be\label{eq:theta-pilot-tuple}([x_1]+\tau_1,\ldots,[x_{\ells}]+\tau_{\ells})\in B^\ells
\ee  which appears in Theorem~\ref{thm:theta-pilot-object-appears}({\bf5}), is the analog, in my theory, of $\Theta$-pilot object considered in \cite[Page 420]{mochizuki-iut3}. 
\subpara The set  of tuples \eqref{eq:theta-pilot-tuple} is also a function of the point of $\sigfjxecpt$ chosen for its calculation. At any rate,  the point of $\sigfjxecpt$ which gives rise to the tuple in \eqref{eq:eq:log-theta-pilot-variation} (i.e. the tuple \eqref{eq:theta-pilot-tuple}) in Theorem~\ref{thm:theta-pilot-object-appears}  can itself be considered to be a variable in \eqref{eq:eq:log-theta-pilot-variation}. 

Thus it makes perfect sense to consider all such tuples so obtained, by  varying the elements of $\sigfjxecpt$. This is done below in \ssep\ref{pa:bound-of-cor312}. \subpara The   quantity 
\be\label{eq:approx-theta-pilot-object}
\frac{1}{\ells}\sum_{j=1}^{\ells}j^2\cdot\log\abs{\xi}_{\C_p}
\ee
appears in the definition of  the ``arithmetic degree of the hull of  $\theta$-pilot object''  in \moccor\ and that degree is a key ingredient in that corollary.
\subpara \emph{Strictly speaking, for Diophantine applications, one should do this for all primes of semi-stable reduction simultaneously.} This can be done with considerably more footwork and will be detailed in \cite{joshi-teich-rosetta}. So Theorem~\ref{thm:theta-pilot-object-appears} should be considered to be a prototype for global constructions.
\subpara Note that the terms on the left hand side of \eqref{eq:log-theta-pilot} are evidently mixed by the actions of  symmetries available in Arithmetic Teichmuller Theory of \cite{joshi-teich}. So one can enlarge the set to incorporate stability with respect to these symmetries. This is done in Definition~\ref{def:basic-theta-values-locus} below.

\nwss

\para\label{pa:bound-of-cor312} 
Now one is ready to define my version $\ttheta$ of the set $\ttheta_{Mochizuki}$ of the theta-values. 
\begin{defn}\label{def:basic-theta-values-locus}
Let $X/E$ be as in \ssep\ref{se:theta-values}. Then the \textit{basic theta-values locus}, denoted by $\ttheta$, is given as follows: let 
\be\ttheta\subset (B)^\ells=\underbrace{B\times B\times \cdots\times B}_{\ells-\text{factors}}\ee 
be the smallest possible subset of $B^\ells$ containing all the tuples \eqref{eq:theta-pilot-tuple}  $$([x_1]+\tau_1,\ldots, [x_\ells]+\tau_\ells)\in B^\ells$$ 
which are lifts $\xi_1$ of the form $[x_j]+\tau_j\in B$ obtained from some point $\sigfjxecpt$ such that $\ttheta$ is both
\benumlab
\item stable under the action of Galois on $B$
\item stable under the action of Frobenius on $B$.
\eenum
\end{defn}
\para Let $B^+\subset B$ be the subring defined in \cite[Chapitre 1, 1.10]{fargues-fontaine}, equipped with the action of Galois and Frobenius $\vphi:B^+\to B^+$.
\blem
Let $B^+\subset B$ be the  subring of $B$ defined in \cite{fargues-fontaine}. Then $$\ttheta\subset (B^+)^\ells\subset B^\ells.$$
\elem
\bp 
Let $z=(z_1,\ldots,z_\ells)\in\ttheta$ corresponding to $(y_1,\ldots,y_\ells)\in\tsigcpt$ with the property that $\eta_{K_j}(z_j)=\xi_1$ i.e. $z_j$ is the lift of $\xi_1\in K_j$ to $B$ (here $K_j$ is the residue field of $y_j$). Let $\rho=1$ and consider the corresponding norm $\abs{-}_\rho=\abs{-}_1$ on $B$. 

Then I claim that for each $1\leq i\leq \ells$ one has $$\abs{z_i}_1\leq 1.$$
Indeed, any $z_i$  considered here is of the form $z_i=[x]+\lambda_i$ for some $\lambda_i\in T_{y_i}\subset B^{\vphi=p}$.  By \cite[Proposition 4.1.3]{fargues-fontaine} one has $$B^{\vphi=p}=(B^+)^{\vphi=p}\subset B^+\subset B$$ and by \cite[Proposition 1.10.7]{fargues-fontaine} one has $$B^+=  \left\{\lambda\in B: \abs{\lambda}_1\leq 1\right\}.$$ 
So $\lambda_i$ satisfies $\abs{\lambda_i}_1\leq 1$ and hence one sees that 
$$\abs{z_i}_1=\abs{[x]+\lambda_i}_1\leq \max\left(\abs{[x]}_1,\abs{\lambda_i}_1\right)\leq \max\left(\abs{[x]}_1,1\right).$$
But by the definition of $\abs{-}_1$ (\cite[D\'efinition 1.4.1]{fargues-fontaine}) one has
$$\abs{[x]}_1 = \abs{x}_F <1 \text{ for any }x\in\fm_F.$$ 
This shows that $\abs{z_i}_1\leq 1$. Since $B^+$ is stable under the action of Galois and the action of Frobenius of $B$, one sees that $\ttheta\subset (B^+)^\ells\subset B^\ells$.
This proves the claim.
\ep
\newcommand{\bcris}{B_{cris}}
\newcommand{\bst}{B_{st}}
\newcommand{\bdr}{B_{dR}}
\para The following is an important consequence of the properties of the ring $B$:
\bpro 
One may naturally view $\ttheta\subset \bcris^\ells \subset \bst^\ells \subset \bdr^\ells$.
\epro
\bp 
This is immediate  using the inclusions $$B^+\subset \bcris\subset \bst\subset \bdr$$
established in \cite{fargues-fontaine} and \cite{fontaine94b}.
\ep

\brem One may further enlarge $\ttheta$  by requiring that $\ttheta$ be stable under the natural action of $${\rm Aut}(G_E\act \sG(\O_{\cpt})=B^{\vphi=p})\isom {\rm Aut}(G_E\act \sG(\O_{\cpt})={(B^+)}^{\vphi=p}).$$ This is one of the  ways in which Mochizuki's Indeterminacy Ind1 \cite{mochizuki-iut3}  arises in my theory (for more details see  \cite{joshi-gconj} and \cite{joshi-teich-rosetta}).
\erem

\para Since $B$ is a Fr\'echet algebra over $\Q_p$ equipped with a family of norms $\left\{\abs{-}_{\rho}: \text{ for }\rho\in [0,1]\subset \R \right\}$, one has a natural way of measuring the ``size of $\ttheta$'' as follows.
\begin{defn} 
For $\rho\in(0,1)\subset\R$, let 
$$|\ttheta|_{B,\rho}=\sup_z\left\{\prod_{j=1}^\ells\abs{z_i}_\rho: z=(z_1,\cdots,z_\ells)\in \ttheta \right\}$$
and let
$$|\ttheta|_B:=\sup_{\rho\in(0,1)}\left\{|\ttheta|_{B,\rho}\right\}=\sup_{\rho\in(0,1)}\sup_z\left\{\prod_{j=1}^\ells\abs{z_i}_\rho: z=(z_1,\cdots,z_\ells)\in \ttheta \right\}.$$
\end{defn}
\para In the parlance of \iut\, the set $\ttheta$ is the \textit{multi-radial representation of theta-values} in the Arithmetic Teichmuller Theory of \cite{joshi-teich,joshi-untilts} and it corresponds to Mochizuki's multi-radial representation of theta-values $\ttheta_{Mochizuki}$ constructed for one fixed $p$-adic field arising as the completion of a number field at one fixed prime. The set $\ttheta_{Mochizuki}$  is the collection of all the $\theta$-pilot object(s) (i.e. analogs of theta-value tuples \eqref{eq:theta-pilot-tuple}) of Mochizuki's Theory, which views such tuples in Galois cohomology (see the next paragraph \ssep\ref{pa:diff-between-two-sets}).  Mochizuki's construction of $\ttheta_{Mochizuki}$ is adelic and hence one necessarily constructs this at each prime--see Remark~\ref{re:global-sit} for more on the number field case. The adelic version of $\ttheta$ is detailed in \cite{joshi-teich-rosetta}.

\subpara\label{pa:diff-between-two-sets} Note however  the set $\ttheta$ constructed above is not the same as the set $\ttheta_{Mochizuki}$ constructed in \iut\ because  Mochizuki does not work with the ring $B$ but works with log-shells tensored with $\Q_p$. As proved in  \cite[Remark 17.5]{joshi-anabelomorphy}, Mochizuki's log-shell tensored with $\Q_p$ can be naturally identified with the Fontaine subspace $H^1_f(G_E,\Q_p(1))\subset H^1(G_E,\Q_p(1))$. By the well-known close relationship (as detailed in \cite[Chapitre 10]{fargues-fontaine}) between these cohomology groups and the cohomology of $G_E$ with coefficients in $B$ or special subspaces of $B$,  my results may be translated quite directly to Mochizuki's version. 

In particular, after the results of this paper,  there is little reason to doubt that $\ttheta_{Mochizuki}$ of \cite{mochizuki-iut3} can be constructed using Mochizuki's prescription but using  Mochizuki's Ansatz established here (as the natural $\Theta_{gau}$-Link) and the theory arithmetic holomorphic structures established in \cite{joshi-teich,joshi-untilts}. 

\subpara\label{pa:tensor-prod-1} Now let me elaborate the significant advantage which one obtains in working with the ring $B$ as opposed to working with Galois cohomology $H^1_f(G_E,\Q_p(1))$ or $H^1(G_E,\Q_p(1))$. The tuples \eqref{eq:theta-pilot-tuple} $([x_1]+\tau_1,\ldots,[x_{\ells}]+\tau_{\ells})\in B^\ells$ by construction. As $B$ is a ring  and so one has a natural multiplication morphism $$B^\ells \to B$$
given by $(\alpha_1,\ldots,\alpha_\ells)\mapsto \alpha_1\cdot\alpha_2\cdots\alpha_\ells$. In particular it makes sense to multiply the individual lifts 
$$([x_1]+\tau_1)\cdot([x_2]+\tau_2)\cdots([x_{\ells}]+\tau_{\ells}) \in B.$$ 
As the multiplication map $B^\ells \to B$ factors naturally as
$$B^\ells\to \underbrace{B\tensor_{\Q_p} B\tensor_{\Q_p}\cdots \tensor_{\Q_p} B }_{\ells-\text{factors}}\to B$$
one may also work with the tensor products of the Fr\'echet algebra $B$ and its variants. This is how tensor product appear naturally in the theory from my point of view. 
In \cite{mochizuki-iut3}, one lacks the ring structure provided here by $B$. Mochizuki forces this product structure by taking tensor products of vector spaces $H^1(G_E,\Q_p(1))$ (see  \cite[\ssep3]{mochizuki-iut3} where this is described as ``tensor-packet structure''). 
\subpara\label{pa:tensor-prod-2} Let me remark that the product of valuations appearing in \eqref{eq:log-theta-pilot-prod-ver} arises from tensor product semi-norms (resp. norms) on tensor products of Banach or Fr\'echet spaces over $p$-adic fields (as in \ssep\ref{pa:tensor-prod-1})  and these tensor semi-norms (resp. norms) satisfy the important property (see \cite[Chap IV]{schneider-book}):
\be\label{eq:cross-norm-property}\abs{v\tensor w}_{V\tensor_{\Q_p}W}=\abs{v}_V\cdot \abs{w}_W\ee  which is referred to as   the \textit{cross-norm property} in the literature on  topological vector spaces over archimedean fields (see \cite{kaniuth-banach-algebras-book}, \cite{grothendieck-memoir} or \cite{diestel-grothendieck-tensor-products-book}).  This is why products of absolute values appearing in \eqref{eq:log-theta-pilot-prod-ver} are quite natural! 

\subpara An important point which emerges from my approach to constructing $\ttheta$ is that as $B$ is a Fr\'echet space, any choice of a norm on $B$ gives us a way of measuring the difference between any two tuples \eqref{eq:theta-pilot-tuple}. This, in some sense, is essential even for \iut\  as its claim rests on being able to (metrically) distinguish between the $\theta$-pilot objects constructed in \cite{mochizuki-iut3}.

\section{A Prototype of Mochizuki's Corollary~3.12}\label{se:proof-of-cor312}
\para Before proceeding to the main result below, let me remark that \moccor\ is dependent on \cite[Theorem 3.11]{mochizuki-iut3} which establishes the existence of multi-radiality i.e. existence of arithmetic holomorphic structures (in the sense of \cite{mochizuki-iut2}) satisfying additional properties. In my case the definition of arithmetic Teichmuller space is based on my notion of arithmetic holomorphic structures (and which provides Mochizuki's version of arithmetic holomorphic structures) \cite{joshi-teich,joshi-untilts} hence multi-radiality is built into the construction of Arithmetic Teichmuller spaces \cite{joshi-teich,joshi-untilts}. It is precisely because of this that one can proceed to the discussion of $\ttheta$ and its properties (such as \moccor) without separately having to establish \cite[Theorem 3.11]{mochizuki-iut3}.

\para 
The following bound is the prototype of \moccor\ for the case of one $p$-adic field.

\bthm\label{th:main} Let $X/E$ be as in \ssep\ref{se:theta-values}. Let $q_E$ be the Tate parameter of $X/E$.
If $\ell$ is a sufficiently large odd prime number, then one has 
\be\label{eq:inequality-of-cor3.12} |\ttheta|_B\geq |q_E^{1/2\ell}|^\ells_{\C_p}.\ee
\ethm
\bp 
Let $(X/E,(E\into K_{y_0}, K_{y_0}^\flat\isom \cpt),*_{K_{y_0}}:\sM(K_{y_0}\to \xan_E))$ be a holomorphoid of $X/E$ with  $y_0\in\syQp$, a closed classical point lying in the fiber over the canonical point of $\sxQp$ i.e. $\pi(y_0)=x_{can}$ and let $q\in K_{y_0}^*$ be the Tate parameter of $\xan_{K_0}$. Then one has $q=\iota_{y_0}(q_E)\in K_{y_0}$; our choice of absolute value on $E$ will be the restriction of $\abs{-}_{K_{y_0}}$ to the copy of $E$ embedded in $K_{y_0}$ by the embedding $\iota_{y_0}$. So $\abs{q_E}_E=\abs{q}_{K_{y_0}}$.
 
Let $t\in \cpt=F$ be such that $([t]-p)\subset W(\O_F)$ is the prime ideal corresponding to  $y_0\in\syQp$.  Since $\tsigcpt$ is stable under Frobenius by Proposition~\ref{pr:frob-stable-cannonical}, the entire Frobenius orbit of $y_0$ is contained in $\tsigcpt$ (by construction $\ttheta$ is  also stable under Frobenius). The choice of $y_0$ is immaterial in the assertion as $\abs{\ttheta_B}_B$ is a supremum. The residue field of $x_{can}$ (and hence $K_{y_0}$) can be identified by \cite[Th\'eor\`eme 10.1.1]{fargues-fontaine} with $\C_p$ with its natural action of $\gal(\bQ_p/\Q_p)$. Normalize the valuation of $K_{y_0}$ so that $\abs{p}_{K_{y_0}}=1$. [For additional properties of valuations of the residue fields of  $\vphi^n(y_0)$ see Theorem~\ref{th:flog-kummer-correspondence} below.]

Let $q$ be the Tate parameter of the holomorphoid $(X/E,E\into K_{y_0}=\C_p,K_{y_0}^\flat\isom \cpt)$. Let $\ell\geq 3$ be a prime, $\ells=\frac{\ell-1}{2}$. Choose a root $t^{1/\ells^2}\in\cpt$. Let  for  $j=1,\cdots,\ells$,
$$a_j=\left(t^{1/\ells^2}\right)^{j^2}\in\cpt,$$
and prime ideals
$$\wp_j=[a_j]-p.$$
Then one has $$([a_1]-p,[a_2]-p,\ldots,[a_\ells]-p)\in\Sigma_F,$$
and let $K_1,\ldots,K_\ells$ be their residue fields. By construction,  $K_\ells$ is the residue field of the canonical point of $\sxQp$ notably the residue field $K_\ells=\C_p$ (as a $\gal(\bQ_p/\Q_p)$-module). 

Now let us compute the relationship between valuations of $K_1,\ldots,K_\ells$. This can be done by Proposition~\ref{pr:lift-vals}. One has for $j=1,\ldots,\ells$:
$$v_{K_j}(p)=j^2\cdot v_{K_1}(p),$$
and hence in particular one has
$$v_{K_\ells}(p)=v_{\C_p}(p)=1=\ells^2 \cdot v_{K_1}(p),$$
and hence one has
$$v_{K_1}(p)=\frac{1}{\ells^2}\cdot v_{\C_p}(p).$$
So one has
\be\label{eq:val-rel-main} v_{K_j}(p)= j^2\cdot v_{K_1}(p)=\frac{j^2}{\ells^2}\cdot v_{\C_p}(p).
\ee

Now let $[x_j]\in W(\O_F)\subset B$ be a Teichmuller lift of $\xi=\xi_1\in K_j$ with $\abs{\xi}_{K_j}=\abs{x_j}_F$ (as in Proposition~\ref{pr:teichmuller-lift-A}, \ref{pr:teichmuller-lift-B}) and Theorem~\ref{thm:theta-pilot-object-appears}.
Then
\be\label{eq:val-rel-main2} 
\sum_{j=1}^\ells\log\abs{[x_j]}_\rho=\sum_{j=1}^\ells \log\abs{\xi}_{K_j},
\ee
and by the established relationship between valuations \eqref{eq:val-rel-main} one has
$$\abs{\xi}_{K_j}=\abs{\xi}_{\C_p}^{\frac{j^2}{\ells^2}}.$$
So one can write the right hand side of \eqref{eq:val-rel-main2} as
$$\sum_{j=1}^\ells\log\abs{[x_j]}_\rho=\sum_{j=1}^\ells \frac{j^2}{\ells^2}\log\abs{\xi}_{\C_p}.$$

This gives 

$$\sum_{j=1}^\ells\log\abs{[x_j]}_\rho=\sum_{j=1}^\ells \frac{j^2}{\ells^2}\log\abs{\xi}_{\C_p}=\frac{\ells(\ells+1)(2\ells+1)}{6\ells^2}\log\abs{\xi}_{\C_p}.$$

Now using $$\abs{\xi}_{\C_p}=\abs{q^{1/2\ell}}_{\C_p},$$ the right hand side of the above equation simplifies to
$$\frac{\ells(\ells+1)(2\ells+1)}{6\ells^2}\cdot\frac{1}{2\ell}\log\abs{q}_{\C_p}.$$
Using $\ells=\frac{\ell-1}{2}$, or equivalently $2\ells+1=\ell$ and canceling factors  this simplifies to
$$\frac{1}{12}\left(1+\frac{1}{\ells}\right)\log(\abs{q}_{\C_p}).$$
If $\ell$ is sufficiently large then $$\frac{1}{12}\left(1+\frac{1}{\ells}\right)<\frac{1}{4}\left(1-\frac{1}{\ell}\right)=\frac{\ells}{2\ell},$$ 
and as $\log(\abs{q}_{\C_p})<0$ one gets
$$\frac{1}{12}\left(1+\frac{1}{\ells}\right) \cdot \log(\abs{q}_{\C_p})> \frac{1}{4}\left(1-\frac{1}{\ell}\right)\cdot \log(\abs{q}_{\C_p})=\log(\abs{q^{1/2\ell}}^\ells_{\C_p}).$$
Thus $\Sigma_F$ contains a point at which $$\sum_{j=1}^\ells\log(\abs{[x_j]}_{\rho})>\log(\abs{q^{1/2\ell}}^\ells_{\C_p}).$$
Since $\log(|\ttheta|_B)$ is the supremum over $\Sigma_F$ (or over a bigger set), one obtains the asserted bound
$$|\ttheta|_B \geq \abs{q^{1/2\ell}}^\ells_{\C_p}.$$ 
This proves the theorem.
\ep

\bcor\label{co:inequality-of-cor3.12-v2} 
In the above notation, if $c\in \R^*$ is a non-zero real number such that 
\be\label{eq:inequality-of-cor3.12-v2}|\ttheta|_B\leq c|q^{1/2\ell}|_{\C_p}^{\ells}\ee
then $$c\geq 1.$$
\ecor 
\bp 
If $c<1$ then this contradicts the lower bound provided by Theorem~\ref{th:main}.
\ep

\para Let me remark that the inequality \eqref{eq:inequality-of-cor3.12} and especially the inequality \eqref{eq:inequality-of-cor3.12-v2} proved in Theorem~\ref{th:main} and Corollary~\ref{co:inequality-of-cor3.12-v2} is a version (in the theory of \cite{joshi-teich}) of the inequality asserted in \cite[Corollary 3.12, Page 598]{mochizuki-iut3}. To see that this is indeed the case, note that, if $0<x<1$ is a real number then $\log(x)<0$ and so, by definition
$$\abs{\log(x)}=-\log(x).$$
As $q$ is the Tate parameter, one has $\abs{q}_{\C_p}<1$ and hence $\abs{q^{1/2\ell}}_{\C_p}<1$. Now suppose that $|\ttheta|_B<1$. For the bound of \eqref{eq:inequality-of-cor3.12} to be useful \emph{at all} this should certainly need to be the case--also see Proposition~\ref{pr:trivial-upper-bound-on-ttheta}. Then $$|\log(|\ttheta_B|)|=-\log(|\ttheta_B|)$$ and similarly $$|\log(|q^{1/2\ell}|_{\C_p}^\ells)|=-\log(|q^{1/2\ell}|_{\C_p}^\ells)$$
and so the inequality \eqref{eq:inequality-of-cor3.12} can also be written as
\be\label{eq:inequality-of-cor3.12-2}
-|\log(|\ttheta|_B)|\geq -|\log(|q^{1/2\ell}|_{\C_p}^\ells)|.
\ee
which is precisely the sort of inequality in \moccor, especially if one writes $$|\ttheta_{Mochizuki}|_B=\sqrt[\ells]{|\ttheta|_B},$$ or equivalently
$$\log(|\ttheta_{Mochizuki}|_B)=\frac{1}{\ells}\log(|\ttheta|_B),$$
then one gets from \eqref{eq:inequality-of-cor3.12-2} that
\be 
-|\log(|\ttheta_{Mochizuki}|_B)|\geq -|\log(|q^{1/2\ell}|_{\C_p})|.
\ee
which has the same shape as the inequality of \moccor. 

Since $\ttheta_{Mochizuki}$ is different from the set $\ttheta$ constructed here (see my discussion \ssep\ref{pa:diff-between-two-sets}) so one cannot assert that \eqref{eq:inequality-of-cor3.12-2} is exactly the same inequality as \moccor. But the inequality stated here is the exact  analog of Mochizuki's inequality in the theory of \cite{joshi-teich,joshi-untilts}.

\brem For  discussions of  \moccor\ from Mochizuki's point of view and his formalism  see \cite{mochizuki-iut3}, \cite{fucheng}, \cite{yamashita}, \cite{dupuy2020statement} (where \moccor\ is stated as Conjecture 1.0.1), \cite{dupuy2020probabilistic} and \cite{dupuy2021kummer}.
\erem

\brem\label{re:global-sit} At this juncture some comments on the number field situation are merited. Suppose $X/L$ is a geometrically connected, smooth hyperbolic curve of type $(g,r)=(1,1)$ over a number field $L$.  Mochizuki's choice of ``initial theta-data'' as in \cite[Definition 3.1]{mochizuki-iut1} fixes another number field $L'$ which is an overfield $L'\supset L$  (the notational translation is $L=F$ and $L'=K$). 

In this setup one wants to bundle together (in a multiplicative sense) the theta-values for all primes $w$ of $L'$ lying over a fixed prime $v$ of $L$ (from my point of view this  precisely amounts to working with $\prod_{w|v}B_{L'_w}$ and hence with tensor products $\bigotimes_{w|v}B_{L'_w}$) (this is a more general version of the discussion of \ssep\ref{pa:tensor-prod-1}). 

There are some technical hurdles which arise because one has to work with tensor products of infinite dimensional F\'rechet algebras as opposed to tensor products of finite dimensional Galois cohomology (which Mochizuki does in \cite{mochizuki-iut3}). 

But at any rate, one arrives at the same mathematical tensor product structure asserted in \cite[\ssep 3]{mochizuki-iut3}--this is detailed in \cite{joshi-teich-rosetta}.  

So at this point,  subject to verification of this paper, \cite{joshi-teich-rosetta} and my other relevant works and modifications to \iut\ suggested therein, there is little reason (\textit{for me}) to doubt the veracity of \moccor.
\erem 

\para The following proposition shows that it is not unreasonable to expect that $|\ttheta|_B$ is bounded from above. For example one has the following upper bound:
\bpro\label{pr:trivial-upper-bound-on-ttheta} 
For the set $\ttheta$ constructed in \ssep\ref{se:construction-ttheta} and for $\rho=1$, one has:
\be 
|\ttheta|_{B,1}\leq 1.
\ee
\epro
\bp This is immediate from the fact that $\ttheta\subset B^+$ and  by \cite[Proposition 1.10.7]{fargues-fontaine} one has $$B^+=  \left\{\lambda\in B: \abs{\lambda}_1\leq 1\right\}.$$ 
This proves the claim.
\ep

\section{Appendix: Mochizuki's $\log$-links and $\flog$-links}\label{se:log-links}
\para I will now demonstrate (in Theorem~\ref{th:log-link}, Theorem~\ref{th:flog-links} and Theorem~\ref{th:flog-kummer-correspondence}) that precise versions of Mochizuki's  $\log$-links, $\flog$-links and Mochizuki's $\flog$-Kummer Indeterminacies are present in the theory of \cite{joshi-teich,joshi-untilts}. Mochizuki's theory of $\log$-links is detailed in \cite{mochizuki-topics3}, \iut--especially \cite[Page 24]{mochizuki-iut3}; other treatments are available--\cite[Section 6.1, Definition 6.3]{fucheng}, \cite{yamashita} and \cite{dupuy2020statement}. 

I will not recall Mochizuki's Theory of $\log$-links since the approach given below (Theorem~\ref{th:log-link}) provides a cleaner demonstration of an  assertion crucial  to \iut\ namely that a  $\log$-link intertwines two possibly distinct (topological) field structures. A $\flog$-link will be a special type of  $\log$-link. 

The main new result here is a precise identification of Mochizuki's $\flog$-link (Theorem~\ref{th:flog-links}) using the Frobenius morphism $\vphi:\syQp\to\syQp$ (this is an adic space over $\Q_p$). This allows one to identify the vertical column of $\flog$-links in \moccor\ with the fiber  of the morphism $\syQp\to \sxQp$ over the canonical point of $\sxQp$.

In particular, Mochizuki's $\flog$-Kummer Indeterminacy Ind3 emerges quite naturally from my point of view (Theorem~\ref{th:flog-links}). As Mochizuki reminds us in \cite{mochizuki-essential-logic}, Indetermincay Ind3 is central to his proof of \moccor.

\para Let $X/E$ be a geometrically connected, smooth, quasi-projective variety over a $p$-adic field $E$. Let $(X/E,(E\into K_y, K_y^\flat\isom \cpt), *_{K_y}:\sM(K_y)\to \xan_{E}))\in \fjxe_{\cpt}$ be a holomorphoid of $X/E$ for some closed classical point $y\in\syQp$.

\para Consider $\hgm$ as a formal (Lubin-Tate) group over $\Z_p$ and let $E=\Q_p$, $F=\cpt$ and let $\sG/\Z_p$ be the Lubin-Tate group with formal logarithm $\sum_n\frac{T^{p^n}}{p^n}$. Let $\sgt$ be the universal cover of $\sG$ given by $\sgt$ as the limit of diagrams (see \cite[Chap. 4]{fargues-fontaine}) $$\sgt=\cdots\mapright{p}\sG\mapright{p} \sG\mapright{p}\sG.$$
\para For $\sG=\hgm$ this can be described quite concretely as follows: $$\hgm(\O_{K_y})=1+\fm_{K_y}$$ and so
$$\widetilde{\hgm(\O_{K_y})}=\left\{(x_n): x_n\in 1+\fm_K,  x_{n+1}^p=x_n  \right\}$$ and by \cite[Proposition 4.5.11]{fargues-fontaine}
$$\widetilde{\hgm(\O_{K_y})} \isom \hgm(\O_{\cpt})=1+\fm_{\cpt}$$
and for simplicity of notation, I will habitually write $$\widetilde{1+\fm_{K_y}}=\widetilde{\hgm(\O_{K_y})}.$$
By \cite{fargues-fontaine}, one has an identification
$$B^{\vphi=p}\isom \sgt(\O_{K_y})=\widetilde{1+\fm_{K_y}} \isom  \hgm(\O_{\cpt})=1+\fm_{\cpt}.$$
\para Let $\eta_{K_y}:B\to K_y$ be the canonical surjection given by the closed classical point $y\in\syQp$.
One has a surjection $$\eta_{K_y}:B^{\vphi=p}\isom \hgm(\O_{\cpt})\to K_y$$ given by the restriction of $\eta_{K_y}$ to $B^{\vphi=p}$. By \cite{fargues-fontaine} this restriction $\eta_{K_y}\big\vert_{B^{\vphi}=p}$  is the $p$-adic logarithm $$B^{\vphi=p}\isom \hgm(\O_{\cpt})=1+\fm_{\cpt} \mapright{\eta_{K_y}=\log} K_y$$ with values in $K_y$   given by $$1+\fm_{\cpt} \ni x\mapsto \log([x])\in K_y,$$
where $[-]:\O_{\cpt}\to W(\O_{\cpt})$ is the Teichmuller lift (\cite[Example 4.4.7]{fargues-fontaine}). 

\para Now suppose one has two holomorphoids $(X/E,(E\into K_{y_1}, K_{y_1}^\flat\isom \cpt), *_{K_{y_1}}:\sM(K_{y_1})\to \xan_{E}))\in \fjxe_{\cpt}$ and $(X/E,(E\into K_{y_2}, K_{y_2}^\flat\isom \cpt), *_{K_{y_2}}:\sM(K_{y_2})\to \xan_{E}))\in \fjxe_{\cpt}$. Then one obtains the diagram of topological vector spaces and continuous linear mappings
$$K_{y_1}\mapleft{\eta_{K_{y_1}}=\log} B^{\vphi=p}\isom \hgm(\O_{\cpt})\isom B^{\vphi=p} \mapright{\eta_{K_{y_2}}=\log} K_{y_2}$$
with the mappings on the extreme left and right being surjections.
\para 
In particular, as $\log$ is not a ring homomorphism and by \cite{kedlaya18} the fields $K_{y_1}$ and $K_{y_2}$ need not be isomorphic, nevertheless the fields $K_{y_1}$ and $K_{y_2}$ are related through the isomorphism of the unit groups \be\label{eq:def-log-link}K_{y_1}\mapleft{\log}\widetilde{\hgm(\O_{K_{y_1}})}\isom \hgm(\O_{\cpt}) \isom B^{\vphi=p} \isom  \hgm(\O_{\cpt})\isom \widetilde{\hgm(\O_{K_{y_2}})}\mapright{\log} K_{y_2}.\ee

\para
In the parlance of \iut\ such a relationship between algebraically closed valued fields containing an isometrically embedded $\Q_p$ is called a \textit{$\log$-link}. \textit{Note however} that Mochizuki does not work with perfectoid fields and works with algebraic closures of $p$-adic fields. Let me provide his version from my point view. Let $\bE_1\subset K_{y_1}$ for the embedding $\iota_{y_1}:E\into K_{y_1}$ provided by the datum of $y_1$ (resp. $\bE_2\subset K_{y_2}$ for the embedding $\iota_{y_2}:E\into K_{y_2}$ provided by the datum of $y_2$) be the algebraic closures of $E$ in $K_{y_1}$ and $K_{y_2}$ (for the respective embeddings). Then Mochizuki's $\log$-link is an abstract isomorphism
$$\O_{\bE_1}^{*pf}\isom \O_{\bE_2}^{*pf}$$ and the surjections
$$\bE_1\mapleft{\log}\O_{\bE_1}^{*pf}\isom \O_{\bE_2}^{*pf}\mapright{\log}\bE_2,$$
where $\O_{\bE_1}^{*pf}=\O_{\bE_1}^*\tensor_\Z\Q$ (resp. $\O_{\bE_2}^{*pf}$) is the perfection of the unit group $\O_{\bE_1}^*$ (resp. $\O_{\bE_1}^*$). In my theory, perfections, in Mochizuki's sense,  are  replaced by $$ \widetilde{\O_{K_{y_1}}^*}=\invlim_{x\mapsto x^p}\O_{K_{y_1}}^* \text{ resp. } \widetilde{\O_{K_{y_2}}^*}=\invlim_{x\mapsto x^p}\O_{\bE_2}^*.$$

\para These constructions and observations prove the following:
\bthm\label{th:log-link} Let 
$(X/E,(E\into K_{y_1}, K_{y_1}^\flat\isom \cpt), *_{K_{y_1}}:\sM(K_{y_1})\to \xan_{E}))\in \fjxe_{\cpt}$ and $(X/E,(E\into K_{y_2}, K_{y_2}^\flat\isom \cpt), *_{K_{y_2}}:\sM(K_{y_2})\to \xan_{E}))\in \fjxe_{\cpt}$ be two holomorphoids of $X/E$. Then one obtains $\log$-linked (in the sense of \eqref{eq:def-log-link}) algebraically closed perfectoid fields $K_{y_1}$ and $K_{y_2}$ and embedded algebraic closures $K_1\supset\bE_1$ (resp. $K_2\supset\bE_2$) of $E$. 
Moreover this gives $\log$-linked isomorphs of the tempered fundamental groups with their natural actions on unit groups
$$\O_{\bE_1}^*\curvearrowleft\pi_1^{temp}(X/E,*_{K_{y_1}}:\sM(K_{y_1})\to \xan_{E}) \isom \pi_1^{temp}(X/E,*_{K_{y_2}}:\sM(K_{y_2})\to \xan_{E})\act \O_{\bE_2}^*.$$
\ethm

\para Now let me come to the construction of $\flog$-links. This is a special case of the above general construction which is central to \moccor. Notably the vertical column of $\flog$-links in \cite[Theorem 3.11]{mochizuki-iut3} and in \moccor\ is a iteration of $\flog$-links. I want to detail how these structures appear from my point of view. 

\para Let $\pi:\syQp\to \sxQp$ be the canonical quotient morphism which identifies the complete curve Fargues-Fontaine curve $\sxQp=\syQp/\vphi^\Z$ as quotient of (adic) curve $\syQp$ by the group generated by Frobenius morphism.  

\para By \cite[Th\'eor\`eme 10.1.1]{fargues-fontaine}, $\sxQp$ has a canonical closed classical point $x_{can}\in\sxQp$ whose residue field may be identified with $\C_p$ equipped with its natural action of $G_{\Q_p}$ (and its open subgroups such as $G_E$). 

\para Let $\{y_n:n\in\Z \}$ be the fiber over $x_{can}$ i.e. $y_n=\vphi^n(y)$ for some $y\in\syQp$ with $\pi(y)=x_{can}$. Notably 
$$y_n=\vphi(y_{n-1})\textit{ for all } n\in\Z.$$ If $\fm_{y_n}\subset B$ is the maximal ideal of $B$ corresponding to $y_n\in\syQp$ then 
$$\vphi(\fm_{y_{n-1}})=\fm_{\vphi(y_{n-1})} =\fm_{y_n}.$$
This is seen as follows. Let $([a]-p)\subset W(\O_{\cpt})$ be the prime ideal corresponding to $y_{n-1}$. Then $\vphi([a]-p)=[\vphi(a)]-p=[a^p]-p$ generates the prime ideal corresponding to $\vphi(y_{n-1})=y_n$. The maximal ideal of $y_{n-1}$ (resp. $y_n$) in $B$ are the closures of these ideals of $W(\O_{ \cpt})$ in $B$.
\para By \cite[Th\'eor\`eme 6.5.2(5)]{fargues-fontaine}, one has an isomorphism 
$$K_{x_{can}}\isom K_{y_n} \textit{ for all }n\in\Z.$$
So one has a natural identification 
$$K_{y_n}\isom \C_p.$$
\para On the other hand, each $y_n$ represents a distinct pair $(K_{y_n},K_{y_n}^\flat\isom \cpt)$ notably the fields $K_{y_n}=\C_p$ is fixed, but the tilting data $(K_{y_n},K_{y_n}^\flat\isom \cpt)$ is not fixed. These observations are summarized in the following:
\bthm\label{th:flog-links} Let $$\left\{(X/E,(E\into K_{y_n}, K_{y_n}^\flat\isom \cpt), *_{K_{y_n}}:\sM(K_{y_n})\to \xan_{E}))\right\}_{n\in\Z}$$ be a collection of holomorphoids of $X/E$ indexed by integers such that the $\{y_n\}_{n\in\Z}$ is the fiber over the canonical point $x_{can}\in\sxQp$. Then 
\benumlab 
\item this collection of holomorphoids is $\flog$-linked in the sense of \eqref{eq:def-log-link}.
\item  Moreover in the  above notation one has a diagram in which all the vertical arrows are continuous surjections given by $p$-adic $\log$ and horizontal arrows are given by the Frobenius morphism:
$$ \xymatrix{
	\cdots \ar[r]^\vphi &  \hgm(\O_{\cpt}) \ar[r]^\vphi\ar[d]^{\log} &  \hgm(\O_{\cpt})  \ar[r]^\vphi\ar[d]^{\log} & \hgm(\O_{\cpt})\ar[r]^\vphi \ar[d]^{\log} & \cdots \\
	\cdots  &  K_{y_{n-1}}\ar[d]^\isom  &  K_{y_{n}}\ar[d]^\isom  & K_{y_{n+1}}\ar[d]^\isom  & \cdots\\
	\cdots  &  K_{x_{can}}=\C_p  &  K_{x_{can}}=\C_p   & K_{x_{can}}=\C_p  & \cdots}$$
\item There is no continuous field isomorphism $\C_p=K_{y_{n-1}}\isom K_{x_{can}}\to K_{y_{n-1}}\isom K_{x_{can}}=\C_p$ which may be inserted into the middle row of the diagram so that each resulting square commutes.
\item In particular, for Mochizuki style $\log$-shells one has $$\log(\widetilde{1+\fm_{y_{n-1}}})\subset p\cdot\log(\widetilde{1+\fm_{y_{n}}}).$$
\eenum
\ethm
\bp 
The first assertion and the second assertion is clear from the preceding discussion. The existence of the diagram is clear from properties of $\syQp$ and its Frobenius morphism established in \cite[Chapitre 6]{fargues-fontaine}. The third assertion is clear as $\C_p$ cannot be equipped with any field homomorphism of the sort $x\mapsto x^p$ which would be forced by commutativity of the diagram.  The last assertion is trivial as $\log([1+x]^p)=p\cdot\log([1+x])$.
\ep

\para\label{pa:mochizuki-quote} In \cite[Remark 1.4.1]{mochizuki-iut3} (also see \cite[Pages 31,32]{mochizuki-iut1}, \cite[Pages 61,62]{mochizuki-gaussian}) Mochizuki describes the significance of his \textit{$\flog$-link} and $\Theta_{gau}$ i.e. a vertical arrow in his $\flog-\Theta$-lattice in the following words:

\textit{``From the point of view of the analogy between the theory of the present
series of papers and p-adic Teichmuller theory [cf. [AbsTopIII, \ssep I5]], the vertical
arrows (i.e. $\flog$-links) of the log-theta-lattice correspond to the Frobenius morphism in positive
characteristic, whereas the horizontal arrows (i.e. $\Theta_{gau}$-link) of the log-theta-lattice correspond
to the `transition from $p^n\Z/p^{n+1}\Z \to p^{n-1}\Z/p^n\Z$', i.e., the mixed characteristic
extension structure of a ring of Witt vectors.''}

\para Mochizuki's version of the diagram in Theorem~\ref{th:flog-links} is \cite[Page 407]{mochizuki-iut3}
$$\xymatrix{
\ar@{.}[r] &\bullet\ar@{-}[r]^\flog\ar[drr] & \bullet\ar@{-}[r]^\flog\ar[dr] & \bullet\ar@{-}[r]^\flog\ar[d] & \bullet\ar@{-}[r]^\flog\ar[dl] &\bullet\ar@{.}[r]\ar[dll] &\\
 & & &  \circ &  &  &
}$$
where each $\displaystyle\bullet$ represents (from my point of view) a  Frobenioid constructed from $$(X/E,(E\into K_{y_n}, K_{y_n}^\flat\isom \cpt), *_{K_{y_n}}:\sM(K_{y_n})\to \xan_{E})).$$ This Frobenioid is constructed by simply  following its geometric construction in \cite{mochizuki-theta}. Similarly  $\displaystyle\circ$ represents a Frobenioid constructed from the category of tempered coverings of $\xan_E$ given by the datum $$(X/E, *_{K_{x_{can}}}:\sM(K_{x_{can}})\to \xan_{E})),$$ using same geometric prescription for the Frobenioid as above.

Each $\displaystyle\circ$ is Mochizuki's ``coric data'' (for \iut) i.e.  identifying each  datum $(X/E, *_{K_{y_n}}:\sM(K_{y_n})\to \xan_{E}))$ (this is the minimal datum required for defining the tempered fundamental group) with the similar datum $(X/E,*_{K_{x_{can}}}:\sM(K_{x_{can}})\to \xan_{E})$. This is done by means of the identification $K_{y_n}\isom K_{x_{can}}$ given by \cite[Th\'eor\`eme 6.5.2]{fargues-fontaine} (this requires forgetting the tilting data provided by each $y_n$). This  provides the identification of each $\displaystyle\bullet$ with $\displaystyle\circ$. 

Notably this identification means the arithmetic holomorphic structures  in the sense of Mochizuki \cite[Example 1.8]{mochizuki-iut2} for $\displaystyle\bullet$ and $\displaystyle\circ$ are isomorphic--more precisely, the \'etale-like pictures are isomorphic but this is not compatible with $\flog$-links (i.e. Frobenius-like pictures are not isomorphic). This is exactly as asserted by \cite[Remark 1.4.2(i), Page 453--454]{mochizuki-iut3}. On the other hand arithmetic holomorphic structures in the sense of \ssep\ref{se:tate-setup} allows one to distinguish between arithmetic holomorphic structures in \eqref{eq:vertical-tower-of-holomorphoids} in my sense (by simply asserting that $y_n\neq y_m$ if $m\neq n$) (as is discussed in \cite{joshi-untilts}, my definition of arithmetic holomorphic structures provides both \'etale and Frobenius-like pictures in the sense of \iut).

\para Theorem~\ref{th:flog-links} demonstrates that once one has the theory of the present series of papers,  one can in fact recognize $\flog$-links and the points  of Mochizuki's Ansatz (\ssep\ref{se:ansatz}) i.e the $\Theta_{gau}$-links of Mochizuki's Theory,   as actually having the  properties mentioned by Mochizuki in \ssep\ref{pa:mochizuki-quote} instead of being analogies!

\para Thus one has obtained a canonical description of Mochizuki's $\flog$-link and the vertical column of $\flog$-links in \moccor.

\bcor 
In the above notation, the collection of holomorphoids $$\left\{(X/E,(E\into K_{y_n}, K_{y_n}^\flat\isom \cpt), *_{K_{y_n}}:\sM(K_{y_n})\to \xan_{E}))\right\}_{n\in\Z}$$
where $\{y_n \}_{n\in\Z}$ is the fiber of the quotient morphism $\pi:\syQp\to\sxQp$  over the canonical point of $x_{can}\in\sxQp$. Then 
\benumlab
\item for any $n\in\Z$ and any choice of $K_{y_{n-1}}$-geometric basepoints $*:\sM(K_{y_{n-1}})\isom\sM(\C_p)\to \xan_{E}$ and $*:\sM(K_{y_{n}})\isom\sM(\C_p)\to \xan_{E}$ one has an isomorphism between tempered fundamental groups
$$\pi_1(X/E,*:\sM(K_{y_{n-1}})\to \xan_{E})\isom \pi_1(X/E,*:\sM(K_{y_{n}})\to \xan_{E})$$
but such an isomorphism is incompatible  with the arithmetic holomorphic structures providing these groups and especially incompatible with the fact that $\vphi(y_{n-1})=y_n$.
\item Moreover diagram given by Theorem~\ref{th:flog-links} can be identified with a vertical column of $\flog$-links in \cite[\ssep 3]{mochizuki-iut3}.
\eenum
\ecor

\para Let me now explain the significance of these results. 
\bthm\label{th:flog-kummer-correspondence} 
Let 
\be\label{eq:vertical-tower-of-holomorphoids}\left\{(X/E,(E\into K_{y_n}, K_{y_n}^\flat\isom \cpt), *_{K_{y_n}}:\sM(K_{y_n})\to \xan_{E}))\right\}_{n\in\Z}\ee 
be the collections of holomorphoids of $X/E$, where  $\{y_n \}_{n\in\Z}$ is the fiber over $x_{can}\in\sxQp$ of the canonical quotient morphism $\pi:\syQp\to\sxQp$. Let $t\in B$ be a function whose divisor in $\syQp$  $$div(t)=\sum_{n\in\Z}(y_n)$$ is the fiber of $\pi:\syQp\to\sxQp$ lying over the canonical point. Then 
\benumlab
\item Each $y_n$ determines a $G_{\Q_p}$-equivariant embedding $\Z_p(1)\to B\subset\bcris$ whose image is $\Z_p\cdot t\subset \bpip\subset B$. 
\item The embeddings of $\Z_p(1)\to B\subset\bcris$ provided by $y_n$ and $y_{n-1}$ differ as $\Z_p$-submodules by  a $p$-multiple by {\bf(1)}.
\item One has $\abs{p}_{K_{y_{n-1}}}=\abs{p}^{1/p}_{K_{y_n}}$.
\item Hence for any fixed $n\in\Z$, and any given $\varepsilon>0$, we can find an $m=m(\varepsilon)\in\Z$ such that
$$\abs{p}_{K_{y_{n-m}}}=\abs{p}^{1/p^m}_{K_{y_n}} > \abs{p}_{K_{y_n}}^{\varepsilon}.$$
\item In particular the valuations of $p$, and hence of the Tate parameters of the holomorphoids $\left\{(X/E,(E\into K_{y_n}, K_{y_n}^\flat\isom \cpt), *_{K_{y_n}}:\sM(K_{y_n})\to \xan_{E}))\right\}_{n\in\Z}$, grow along the fiber $\{y_n\}_{n\in \Z}$.
\eenum
\ethm
\brem Before providing a proof, it will be worth understanding what is being claimed here and why it is relevant to the Diophantine calculations of \iut.
\benumlab
\item Property {\bf(3)} is mentioned as an analogy in \cite[Remark 3.6.2]{mochizuki-topics3}.
\item On $\syQp$ there is no canonical point analogous to the canonical point  $x_{can}\in\sxQp$ and so one must work with the entire fiber $\{ y_n\}_{n\in\Z}$. This means there is no canonical choice of the homomorphism $\Z_p(1)\to B\subset \bcris$.  \textit{This is Mochizuki's $\flog$-Kummer indeterminacy}. Explicitly $(X/E, (E\into K_{y_{n-1}}),*:\sM(K_{y_{n-1}})\to \xan_E)$ provides a Kummer theory (i.e. an embedding $\Z_p(1)\to B\subset \bcris$) which is distinct from that provided by $(X/E, (E\into K_{y_{n}}),*:\sM(K_{y_{n}})\to \xan_E)$ by a $p$-multiple (for a similar assertion \cite[Example 2.6.1(b)]{mochizuki-gaussian}). As Mochizuki reminded us, through \cite[\ssep 3.11]{mochizuki-essential-logic}, that $\flog$-Kummer Indeterminacy is fundamental to his proof of \cite{mochizuki-iut4}.  So it is important to exhibit how this appears from my point of view.
\item The last assertion of Theorem~\ref{th:flog-kummer-correspondence} is analogous to growth of radii of convergence of solutions of $p$-adic differential equations under Frobenius. 
\eenum 
\erem

\brem  
Let me complement the above remark with the example of a Tate elliptic curve $X/E$. Let $\{y_n \}_{n\in\Z}$ be as above and for simplicity of notation let me write $K_n=K_{y_n}\ (\isom \C_p)$ for the residue field of $y_n$ and $\eta_n$ for $\eta_{K_n}:W(\O_{\cpt})\to \O_{K_n}$ for the canonical surjection.  Let $\xi\in\O_{K_n}$ with $\abs{\xi}_{K_n}<1$. Let $[z_n]\in W(\O_{\cpt})$ be a Teichmuller lift of $0\neq\xi\in \fm_{K_n}\subset \O_{K_n}$ to $W(\O_{ \cpt})$. This situation is summarized in the following diagram of ring homomorphisms
$$\xymatrix{
	&\ar[dl]_{\eta_m} W(\O_{ \cpt}) \ar[dr]^{\eta_{n}} & \\
	\O_{K_{m}} && \O_{K_{n}}
}.$$
For each $m\in \Z$, let $$\xi_m=\eta_{m}([z_n])\in \O_{K_{m}}.$$ Then for each $m\in \Z$, 
$$\xi_m\neq0 \text{ in }\O_{K_m}.$$ This is immediately clear from  \cite[Lemma 2.2.13]{fargues-fontaine} and by Theorem~\ref{th:flog-kummer-correspondence} one has
$\abs{\xi_m}_{K_m}\neq \abs{\xi}_{K_n}$.
Now let $[z_n]$ resp. $[z_m]$ be the lifts of the Teichmuller lifts Tate parameter $q_m$ of $X/K_m$ resp. $q_n$ of $X/K_n$. Then $\eta_n([z_m])\neq 0$ and $\eta_m([z_n])\neq0$.  In particular one can consider the set  $\{\eta_n([z_m])\}_{m\in\Z}\subset K_n$. If $A\subset \fm_{\O_{K_{n}}}$ is a compact, bounded subset of $K_n$ then 
$$\{\eta_n([z_m])\}_{m\in\Z}\cap A=\{\eta_n([z_m]): m\gg-\infty\}.$$
\erem

\brem Let me remark that for applications to Diophantine geometry, it is not enough to work with Arithmetic Teichmuller spaces $\fjxe$ viewed over $\sxQp$ and estimate local height contributions using the canonical point $x_{can}$ and its residue field as there is no canonical valuation on the residue field of $x_{can}$. 

On the other hand each $y_n$ comes with a canonical valuation on its residue field $K_{y_n}\isom \C_p$ but the valuations provided by $y_{n-1}$ and $y_n$ are not the same on $\Q_p$ as proved in Theorem~\ref{th:flog-kummer-correspondence}{\bf(3)}.  

As the theorem shows, working with $\syQp$ reveals another feature of the theory, namely Kummer Theory i.e. the galois equivariant homomorphism $\Z_p(1)\to B\subset \bcris$ acquires dependence on $y_n$. 

This means that a Galois cohomology class for one such Kummer homomorphisms is a $p^m$-multiple of the similar class for another such Kummer theory for some $m\in\Z$.  This is a subtle point in the theory and in some sense Mochizuki appears to have been the first to recognize this delicate issue from his anabelian point of view and that it implies that valuations of Tate parameters rise (in a tower of $\flog$-links in \cite{mochizuki-iut3}). Mochizuki's way of understanding this phenomenon is through his cyclotomic rigidity isomorphism. His discussion of this may be gleaned from \cite[Page 20, Examle 2.6.1]{mochizuki-gaussian}.  

At any rate, this means, for Diophantine applications, one must take all the holomorphoids \eqref{eq:vertical-tower-of-holomorphoids} (as above), together with the valuations, and the Kummer Theories provided by each $y_n$ together; and when one does this, the valuations of the Tate parameters grow naturally as is established in above!

The growth of valuations established above is crucial in understanding how one can hope to arrive at upper bounds in the Diophantine contexts by using this idea.
\erem
\bp
The first and the second  assertions are obtained by combining the conclusions of  \cite[Proposition 6.4.1]{fargues-fontaine} and the properties of the Lubin-Tate period mapping \cite[Remarque 6.3.2]{fargues-fontaine} and \cite[Proposition 10.1.1]{fargues-fontaine}. 

Let me prove the third assertion. Suppose the principal prime ideal of $y_n$ in $W(\O_{\cpt})$ is given by $([a]-p)\subset W(\O_{\cpt})$, then  the ideal $([\vphi^{-1}(a)]-p)\subset W(\O_{\cpt})$ is the principal prime ideal corresponding to $y_{n-1}\in\syQp$.

By \cite[Proposition 2.2.17]{fargues-fontaine} one sees that $$\abs{p}_{K_{y_{n-1}}}=\abs{\vphi^{-1}(a)}_{\cpt}=\abs{a}_{\cpt}^{1/p}$$ and $$\abs{p}_{K_{y_{n}}}=\abs{a}_{\cpt}=\abs{a}_{\cpt}$$
or to put it differently 
$$\abs{p}_{K_{y_{n-1}}}=\abs{p}_{K_{y_{n}}}^{1/p}.$$ This proves {\bf(3)}. Finally note that
as $0<\abs{a}_{\cpt}<1$, so $$\abs{p}_{K_{y_{n-1}}}=\abs{a}_{\cpt}^{1/p}>\abs{a}_{\cpt}=\abs{p}_{K_{y_n}}.$$
Now clearly {\bf(3)}$\implies${\bf(4)}. The last assertion is clear from {\bf(4)}.
This completes the proof.
\ep

\para Note that in \iut--especially in \cite{mochizuki-iut3}, Mochizuki works with $$\O_{\bQ_p}^{\times\mu}=\O_{\bQ_p}^*/\mu(\bQ_p)=(1+\fm_{ \bQ_p})/\mu_p(\bQ_p)$$
where $\mu(\bQ_p)$ (resp. $\mu_p(\bQ_p)$) is the subgroup of $\O_{\bQ_p}^*$ of roots of unity  (resp. subgroup of the $p$-power roots of unity) in $\bQ_p$, and especially with the exact sequence
$$1\to \mu_p(\bQ_p)\to 1+\fm_{ \bQ_p} \to \O_{\bQ_p}^{\times\mu}\to 1,$$
and considers ``different Kummer Theories'' (see \cite[Page 39--42]{mochizuki-gaussian}) i.e. different isomorphs of this sequence (and the associated Galois cohomology) obtained via Anabelian Reconstruction Algorithms of \topics, \iut,
while here and in \cite{joshi-teich}, following \cite{fargues-fontaine}, I work with the exact sequence (dependent on $K$)
$$1\to T_p(\hgm)(\O_{K}) \tensor_{\Z_p}\Q_p \to \widetilde{1+\fm_{K}}=\widetilde{\hgm(\O_K)} \to K\to 0,$$
where $T_p(\hgm)(\O_K)$ is the ($p$-adic) Tate-module of $\hgm$ computed using $\O_K$ and $\hgm(\O_K)$ and allow $K$ to vary among algebraically closed perfectoid fields in general. 

Now let $$\bQ_{p;K}\subset K$$ be the algebraic closure of $\Q_p\subset K$ with the valuation induced from $K$. Then one has a Kummer sequence indexed by the algebraically closed perfectoid field $K$:
\be 1\to \mu_p(\bQ_{p;K})\to 1+\fm_{ \bQ_{p;K}} \to \O_{\bQ_{p;K}}^{\times\mu}\to 1.\ee
Especially if $K_1,K_2$ are not topologically isomorphic, then one obtains two distinguishable Kummer sequences of the above type. On the other hand, Theorem~\ref{th:flog-kummer-correspondence} demonstrates that even in the situations where the fields $K_1\isom K_2$, one has distinct Kummer theories in the sense of Theorem~\ref{th:flog-kummer-correspondence} and this distinction impacts geometric data such as the respective Tate-parameters (and their respective valuations) in a rather subtle way.

\bibliography{../../master/masterofallbibs.bib}

\Address
\end{document}